\documentclass{amsart}
\usepackage{graphicx,amssymb,amsmath,epsfig}
\usepackage[all]{xy}
\usepackage[T1]{fontenc}
\usepackage{lmodern}

\usepackage{amssymb}
\usepackage{graphicx}
\usepackage{wrapfig}
\usepackage{amsmath}
\usepackage{amsmath,amssymb}

\usepackage{graphicx}

\newtheorem{theorem}{Theorem}[section]
\newtheorem{proposition}[theorem]{Proposition}
\newtheorem{lemma}[theorem]{Lemma}
\newtheorem{corollary}[theorem]{Corollary}
\newtheorem{remark}[theorem]{Remark}
\newtheorem{example}[theorem]{Example}

\newtheorem{definition}[theorem]{Definition}
\newtheorem{corollaries}[theorem]{Corollaries}
\newcommand{\bth}{\begin{theorem}}
\newcommand{\bpr}{\begin{proposition}}
\newcommand{\epr}{\end{proposition}}
\newcommand{\bco}{\begin{corollary}}
\newcommand{\eco}{\end{corollary}}
\newcommand{\ble}{\begin{lemma}}
\newcommand{\ele}{\end{lemma}}
\newcommand{\bde}{\begin{definition}\rm}
\newcommand{\ede}{\end{definition}\rm}
\newcommand{\bre}{\begin{remark}\rm}
\newcommand{\ere}{\end{remark}}
\newcommand{\bex}{\begin{example}\rm}
\newcommand{\eex}{\end{example}}
\newcommand{\bcors}{\begin{corollaries}\rm}
\newcommand{\ecors}{\end{corollaries}}

\def\ra#1{\hbox to #1pc{\rightarrowfill}}
\def\fract#1#2{\raise3pt\hbox{$ #1 \atop #2 $}}
\def\lrar{{\ra 2}}

\def\sp#1{\hbox{SP}^{#1}}

\def\bsp#1{\overline{\hbox{SP}}^{#1}}
\def\sn{{\mathfrak S_n}}
\def\bbz{{\mathbb Z}}

\def\bbp{{\mathbb P}}
\def\bbr{{\mathbb R}}

\def\bbc{{\mathbb C}}
\def\bbr{{\mathbb R}}

\begin{document}

\title{Homotopy Groups of Diagonal Complements}
\author{Sadok Kallel and Ines Saihi}
\address{Laboratoire Painlev\'e, USTL (France)
and American University of Sharjah (UAE).}
\email{sadok.kallel@math.univ-lille1.fr}

\address{Ecole sup\'erieure des sciences et techniques de Tunis (Tunisia)}
\email{ines.saihi@esstt.rnu.tn}

\keywords{diagonal arrangements, homotopy groups, configuration space, colimit diagram}
\subjclass[2010]{Primary 55Q52; Secondary 55P10,57Q55}

\maketitle

\centerline{In memory of Abbas Bahri so greatly missed}

\begin{abstract}
For $X$ a connected finite simplicial complex we consider $\Delta^d(X,n)$ the space of configurations of $n$ ordered points of $X$ such that no $d+1$ of them are equal, and $B^d(X,n)$ the analogous space of configurations of unordered points. These reduce to the standard configuration spaces of distinct points when $d=1$. We describe the homotopy groups of $\Delta^d(X,n)$ (resp. $B^d(X,n)$) in terms of the homotopy (resp. homology) groups of $X$ through a range which is generally sharp. It is noteworthy that the fundamental group of the configuration space $B^d(X,n)$ abelianizes as soon as we allow points to collide (i.e. $d\geq 2$).
\end{abstract}



\section{Introduction}

Let $X$ be a  topological space and $\Delta_{d+1}(X,n)\subset X^n$ the union of the $(d+1)$-st diagonal arrangement in $X^n$; that is
$$\Delta_{d+1}(X,n) = \{(x_1,\ldots, x_n)\in X^n\ |\
x_{i_0}=x_{i_1}=\cdots =x_{i_d}\ \hbox{for some sequence}\ 1\leq i_0<\cdots < i_d\leq n\}.$$
Its complement in $X^n$ is the ``configuration space of no $d+1$ equal points in $X$'' which is written
$$\Delta^d(X,n) = X^n-\Delta_{d+1}(X,n)\ .$$
This is the space of ordered tuples of $n$ points in $X$ with the multiplicity of each entry in the tuple \textit{at most} $d$ (hence the notation $\Delta^d$ as opposed to $\Delta_d$ for \textit{at least} $d$). It is useful to think of these tuples as ``configurations" of $n$ ordered points in $X$ with the property that $d$ of the points can collide but not $d+1$.
The symmetric group $\mathfrak S_n$ acts on $\Delta^d(X,n)$, and the quotient is denoted by $B^d(X,n)$. We have increasing filtrations
\begin{equation}\label{filtration}
F(X,n) := \Delta^1(X,n)\subset \Delta^2(X,n)\subset\cdots\subset
\Delta^n(X,n)= X^n ,
\end{equation}
$$B(X,n) := B^1(X,n)\subset B^2(X,n)\subset\cdots\subset
B^n(X,n)= \sp{n}X,$$
with $\sp{n}X:=X^n/_{\sn}$ being the $n$-th symmetric product.
Here we have written $F(X,n)$ and $B(X,n)$ for the standard configuration
spaces of ordered (resp. unordered) pairwise distinct points of cardinality $n$. Various other notations for $F(X,n)$ in the literature include $C_n(X)$, $Conf_n(X),\cdots$, while $B(X,n)$ is sometimes written Braid$(X,n)$ in the geometric topology literature; reminiscent of the fact that its fundamental group is the so-called $n$-th braid group of $X$.

In some exceptional cases, the spaces $\Delta^d(X,n)$ and $B^d(X,n)$ can be empty (if for example $X$ is a point and $d<n$), but otherwise they have a rich and interesting geometry (see \cite{kt}). An early appearance of $\Delta^d(X,n)$ is in \cite{cl} in connection with Borsuk-Ulam type results while more recent applications to the \textit{colored Tverberg Theorem for manifolds} appear in \cite{blago}. In the case $V$ is a vector space, the spaces $\Delta^d(V,n)$ are subspace complements dubbed ``non-$(d+1)$-equal arrangements'' in \cite{bw} and their homology is made explicit in \cite{dt} as an algebra over the little disks operad, with interesting applications to the spaces of non-$d$-equal immersions. In the case $X=\bbc$, the spaces $B^d(\bbc,n)$ are intimately related to spaces of based holomorphic maps from the Riemann sphere into complex projective space $\bbp^d$ (see \cite{guest,kallel}). In all cases, these spaces seem to have been studied so far exclusively for when $X$ is a manifold. One of our objectives in this paper is to give some sharp results on the homology and homotopy groups of the non-d equal configurations of $X$ when $X$ is a more general polyhedral space.

Throughout this paper a space $X$ means a finite simplicial complex; that is the realization of a finite abstract simplicial complex. Unless specified, all spaces are connected. 

\bth\label{main1}
Let $X$ be a connected finite simplicial complex that is not a point, $d\geq 2, n\geq 2$. Then
$$\pi_i(B^d(X,n))\cong \pi_i(\sp{n}(X))\ \ \hbox{for}\ \ \ 0\leq i\leq 2d-2\ .$$
In particular $\pi_1(B^d(X,n))\cong H_1(X;\bbz )$ when $d\geq 2$, $n\geq 2$.\\
Moreover if $X$ is simply connected, $2\leq d\leq n$, then
$$\pi_i(B^d(X,n))\cong \tilde H_i(X;\bbz )\ \ \hbox{for}\ \ \ 0\leq i\leq 2d-2\ .$$
where $\tilde H(-;\bbz)$ is reduced integral homology.
\end{theorem}

The bound $2d-2$ in the theorem is sharp as is illustrated by the case $X$ a Euclidean space (\S\ref{main1b}).
Note that the special case of the fundamental group says that ``allowing a single collision is enough to abelianize the fundamental group''. This can be expected since collisions kill the \textit{braiding} (\S\ref{fundamental}).

The homotopy groups of $\Delta^d(X,n)$ turn out to depend on local connectivity properties of the space. We say $X$ has \textit{local homotopical dimension} $r$ if for any $x\in X$ and any neighborhood $U$ of $x$, there is an open neighborhood $V\subset U$ of $x$ such that $V-\{x\}$ is $r$-connected (Definition \ref{defs}).

\bth\label{main2}
Let $X$ be a locally finite simplicial complex with local homotopical dimension $r\geq 0\ ,\ d\geq 1$. Then $$\pi_i(\Delta^d(X,n))\cong \pi_i(X)^n\ \ \ \hbox{for}\ \ i\leq rd + 2d-2\ .$$
\end{theorem}

\bre If $d\geq n$, both spaces are equal $\Delta^d(X,n)=X^n$ and all homotopy groups agree. When $d<n$ this bound is in general optimal as can be seen in the case of manifolds. For example $\bbr^2$ has $0$ local homotopical dimension and $\Delta^d(\bbr^2,d+1)\simeq S^{2d-1}$ is $2d-2$-connected precisely.\ere

\bre For a polyhedral pair $(X,Y)$, the ``homotopical depth" of $Y$ in $X$ is set to be $n$ if the pair $(X,X\setminus Y)$ is $n$-connected \cite{eyral}. Theorem \ref{main2} is saying that the homotopical depth of the diagonal arrangement $\Delta_{d+1}(X,n)$ in $X^n$ is at least $rd + 2d-2$. This seems to be the first complete such calculation for this kind of arrangements of subspaces.
\ere

To prove both theorems we use a localization principle for homotopy groups (Theorem \ref{localizing}) relating the local connectivities of pairs $(V,V\setminus Y)$ to the global connectivity of $(X,X \setminus Y)$ for closed $Y\subset X$ and $V$ local neighborhoods in a cover. In both cases the proof reduces to studying the case of $V$ being the union of various simplices joining along a simplex. For Theorem \ref{main2}, the argument amounts to giving a ``homotopical decomposition" of $\Delta^d(V,n)$ when $V$ is such a union. We recall that by a homotopical decomposition of a space $X$ we mean a \textit{diagram} $\mathcal D: I\lrar Top$; i.e. a functor from a small category $I$ to the category of topological spaces and continuous maps, so that the map $\hbox{hocolim}_I\mathcal D \lrar\hbox{colim}_I\mathcal D\cong X$ is a weak equivalence (see \S\ref{proofthm2}).
 Our decomposition extends similar results in \cite{sun}. Since we are able to control the connectivity of each space making up the diagram, we are able to derive our bound.

Theorem \ref{main1} on the other hand relies on a different argument. First we treat the case of a manifold based on the idea of scanning maps. The general case appeals to a theorem of Smale \cite{smale} relating the connectivity of a map to that of its preimages. Since Smale's theorem works for proper maps, a technical issue we have to deal with is the construction in \S\ref{compactification}
of a $\sn$-equivariant simplicial complex which is a deformation retract of $\Delta^d(X,n)$ for $X$ again a finite complex. As pointed out by the referee, similar techniques are
in (\cite{matroids} chapter 4) and have been applied to hyperplane arrangements in \cite{blago2} for example (see references therein). Section \S\ref{compactification} is of independent interest and has  relevance to more recent constructions of CW-retracts for configuration spaces \cite{jesusdai}.

The first section of the paper discusses motivational examples and general connectivity results. The second section discusses the special case of graphs. Proposition \ref{algorithm} gives a simplified and then expanded version of a useful theorem of Morton, which is used to give an amusing description of the homotopy type of the configuration space of two points on a wedge of circles (Proposition \ref{wedgecircle}).

\vskip 5pt
{\sc Acknowledgment}: The first part of this work was conducted at the University of Lille 1 under a ``BQR" grant. The Mediterranean Institute for the Mathematical Sciences (MIMS) has made resources available during the completion of this work. We are grateful to Faten Labassi for pointing us to Munkres' book and Proposition \ref{munkres}. Finally we thank Paolo Salvatore for his insight on Lemma \ref{null}.

\section{Preliminaries}\label{examples}

We start with some classical examples of diagonal arrangements and their complements. The extreme cases $d=1$ and $d=n-1$ are most encountered in the literature. The case $\Delta^1(X,n)=F(X,n)$ corresponds to the configuration space of pairwise distinct points
$$F(X,n) = \{(x_1,\ldots, x_n)\in X^n\ |\ x_i\neq x_j\ \hbox{for}\ \ i\neq j\}\ .$$
The action of $\sn$ on $F(X,n)$ is free and we have a regular covering $F(X,n)\lrar B(X,n)$.
If $X$ is a manifold of $\dim X >2$, then $\pi_1(F(X,n))\cong\pi_1(X^n)$ by codimension argument (see Proposition \ref{manifoldcase}),
while $\pi_1(B(X,n))$ is a wreath product $\pi_1(X)\wr\mathfrak S_n$  (this is standard but a leisurely exposition is in \cite{imbo}).

\bex\label{casestudy}
When $d=n-1$, $B^{n-1}(X,n)$ is the complement in $\sp{n}(X)$ of the diagonal embedding $\Delta : X\hookrightarrow\sp{n}X$, $x\mapsto [x,\ldots ,x]$. When $X=\bbc$, the elementary symmetric functions give a diffeomorphism $\sp{n}(\bbc )\cong\bbc^n$ and the image of $\Delta (\bbc)$ corresponds under this diffeomorphism to the rational normal curve $V$ diffeomorphic to the Veronese embedding $x\mapsto (x,x^2,\ldots, x^n)$. One can check that
$$B^{n-1}(\bbc, n) \cong \sp{n}(\bbc ) -V\simeq S^{2n-3}\ .$$
A short proof of this equivalence is given in (\cite{guest}, Lemma 2.7), while another quick argument would be to use simple connectivity of $B^{n-1}(\bbc ,n)$ and Alexander duality.

In general for $\bbr^k$, $k\geq 2$,
$\Delta^{n-1}(\bbr^k,n)=(\bbr^k)^n-\Delta$ is the complement of the thin diagonal and this deforms onto the orthogonal complement of the diagonal $\Delta = \{(x,\ldots, x)\}$ minus the origin so that
 $\Delta^{n-1}(\bbr^k,n)$ is up to homotopy the unit sphere
$S^{nk-k-1}$ in $\{(x_1,\ldots, x_n)\in (\bbr^k)^n\ |\ \sum x_i = 0\}=\Delta^\perp$.
 This deformation can be made equivariant with respect to the permutation action of $\sn$
 so that the $\sn$-quotient is $B^{n-1}(\bbr^k,n)$. We show below that this space
 is simply connected as soon as $n\geq 3$ (in fact it is $2n-4$-connected; Lemma \ref{bn-1rk}).
\eex

\ble\label{trivial}
If $S$ is the unit sphere in $H = \{(v_1,\ldots, v_n)\in (\bbr^k)^n,
\sum v_i = 0\}$, and if $\sn$ acts on $H$, and hence on $S$, by permutation of coordinates, then the quotient $Q_{n,k}:=S/\sn$ is simply connected whenever $nk-k-1\geq 2$.
\ele

\begin{proof}
We use the following useful main result of Armstrong \cite{armstrong}:
 let $G$ be a discontinuous group of homeomorphisms of a path connected, simply connected, locally compact metric space $X$, and let $H$ be the normal subgroup of $G$ generated by those elements which have fixed points, then the fundamental group of the orbit space $X/G$ is isomorphic to the factor group $G/H$. We apply this result to $G=\sn$ and $X=S$ which is simply connected. The point is that when $n\geq 3$, the fixed points of the permutation action are of the form $(v_1,\ldots, v_n)$ with $v_i=v_j$ for some $i<j$, which means that all transpositions are in $H$ and hence $G=H$.
\end{proof}

The argument of Armstong used in the proof of Lemma \ref{trivial} implies that if $\Delta^d(X,n)$ is simply connected, then $\pi_1(B^d(X,n))$ is the quotient of $\mathfrak S_n$ by the normal subgroup generated by elements having fixed points and this subgroup is the entire group if $d\geq 2$. This establishes a useful conclusion.

\bco\label{armstrong} If $\Delta^d(X,n)$ is simply connected, then so is $B^d(X,n)$ if $d\geq 2$.
\eco

The following result, valid for smooth (i.e. $C^\infty$) manifolds, is a special case of Theorem \ref{main1}.

\bpr\label{pi1manifold} When $X=M$ is a closed smooth manifold, $\dim m\geq 2$ and $n\geq 3$, then
$\pi_1(B^{n-1}(M,n))$ is isomorphic to $H_1(M;\bbz )$.
\epr

\begin{proof} A tubular neighborhood of the diagonal copy of $M$ in $\sp{n}M$ can be identified with the total space of the following subbundle. Let $TM^{\oplus n}$ be the $n$-fold Whitney sum of the tangent bundle $TM$ of $M$, $\dim M=m$, and let $\eta$ be the subbundle with fiber $H = \{(v_1,\ldots, v_n)\ |\ \sum v_i = 0\}$. The total space of this subbundle is homeomorphic to a neighborhood of diagonal $M$ in $M^{\times n}$. Now $\sn$ acts on this bundle fiberwise (linearly on each fiber) and the fiberwise quotient $\zeta$ has fiber $H/\sn$ which can be identified with the cone
$c(S^{(n-1)m-1}/\sn)$, where $\dim M = m$ and $S^{(n-1)m-1}$ is the unit
sphere in $H$. According to (\cite{kt}, Proposition 4.1) and for $M$ a smooth closed manifold, a neighborhood deformation retract $V$ of the diagonal $M$ in $\sp{n}M$ is homeomorphic to the total space of $\zeta$. The fiberwise apexes of the fiberwise cone give the ``zero-section'' of this bundle. The complement of this section is $S(M)$ which is up to fiberwise equivalence a bundle over $M$ with fiber $S^{(n-1)m-1}/\sn$. By construction we have the homotopy pushout
$$\xymatrix{
S(M)\ar[r]^\pi\ar[d]&M\ar[d]^\Delta\\
B^{n-1}(M,n)\ar[r]^\iota&\sp{n}M}
$$
If $n=2, m\geq 2$, $S(M)$ is the projectivized tangent bundle with fiber $S^{m-1}/\bbz_2 = \bbr P^{m-1}$. When $n\geq 3$ and $m\geq 2$, $S(M)$ has simply connected fiber (Lemma \ref{trivial}) so that $\pi$ induces an isomorphism on fundamental groups and by the Van-Kampen theorem, $\iota$ induces an isomorphism on $\pi_1$ as well; i.e. $\pi_1(B^{n-1}(M,n))\cong \pi_1(\sp{n}M)\cong H_1(M;\bbz )$ for $n\geq 3$.
\end{proof}

To complete this section, we state a well-known result which later will be seen as a special manifestation of the ``localization principle'' (\S\ref{main1b}).

\bpr\label{manifoldcase} If $S=\bigcup S_j$ is a finite union of submanifolds of a smooth manifold $M$, closed with real codimension $d\geq 2$, then the inclusion $M-S\hookrightarrow M$ induces an isomorphism on homotopy groups $\pi_i$ for $0\leq i\leq d-2$, and an epimorphism on $\pi_{d-1}$.  \epr

A proof of the above proposition, using standard transversality arguments, can be found for example
in (\cite{helmke}, Lemma 5.3). This proposition is not true if the ``ambient space'' is not a
manifold. For example $B(\bbr^m,2)$ is the complement of the diagonal in $\sp{2}(\bbr^m)$ and we have the homotopy equivalence $B(\bbr^m,2)\simeq\bbr P^{m-1}$ so that $\pi_1(B(\bbr^{m},2))\cong\bbz_2$ no matter the codimension of the diagonal $m\geq 3$.

As a consequence we have the following precursor of Theorem \ref{main1}

\bco\label{n2d2} If $X$ is a topological surface and $d\geq 2$, then
 $\pi_1(B^d(X,n))\cong H_1(X,\bbz )$.
 \eco

\begin{proof} The real plane $\bbr^2$ has the special property that $\sp{n}(\bbr^2)$ is diffeomorphic to $\bbr^{2n}$. This implies right away that
when $S$ is a topological surface, $\sp{n}(S)$ is a manifold of dimension $2n$,
and that $B_{d+1}(X,n)$ is the union of submanifolds of dimension at most $2(n-d) = 2n-2d$. This means that $B^d(S,n)=\sp{n}(S)-B_{d+1}(S,n)$ is the complement of a finite union of submanifolds of codimension at least $2d>2$. By Proposition \ref{manifoldcase}, $\pi_1(B^{d}(S,n))\cong\pi_1(\sp{n}S)$ and this is again $H_1(S,\bbz)$ for $n>1$.
\end{proof}


\section{The Case of the Circle} \label{circle}

Write $S^1 = \{z\in\bbc , |z|=1\}$ the unit circle in the complex plane. There is a map $B^d(S^1,n)\lrar S^1$ which multiplies the points of a configuration in $S^1$. This map is well defined since $S^1$ is abelian. This map turns out to have contractible fibers so
that in particular $B^d(S^1,n)\simeq S^1$ (see Proposition \ref{bds1}).

Let $\Delta_{n-1}= \{(s_1,\ldots, s_n), 0\leq s_i\leq 1,\sum s_i=1\}$ be the $n-1$ dimensional simplex and write $\Delta_{n-1}(d)$ the partial compactification of the open simplex $\mathring{\Delta}_{n-1}$ where we allow at most $d$ \textit{consecutive} $s_i$'s to be zero (using cyclic ordering, $s_n$ and $s_1$ are consecutive to each other). In particular
$\Delta_n(1) = \mathring{\Delta}_{n-1}$. We will write $\bbz_n$ for the cyclic group of order $n$.
Using a similar action as in (\cite{bm}, p.407) we have:

\bpr\label{algorithm} Let $\bbz_n$ with multiplicative generator $\tau$ act on $S^1\times\Delta_{n-1}(d)$ via
$$\tau (e^{i\theta}, s_1,\ldots, s_n) = (e^{i\theta + i2\pi s_1}, s_2,\ldots, s_n,s_1)\ .$$
Then the quotient by the action; written $S^1\ltimes_{\bbz_n}\Delta_{n-1}(d)$,  is homeomorphic to $B^d(S^1,n)$.
When $d=1$, there is a $\sn$-equivariant homeomorphism
$$F(S^1,n)\cong (S^1\times \mathring{\Delta}_{n-1})\times_{\bbz_n}\sn\ .$$
\epr

\begin{proof}  The cyclic group appears for a simple reason: any configuration $(x_1,\ldots, x_n)$ can be brought into a unique counterclockwise configuration up to cyclic permutation.
More precisely let $(x_1,\ldots, x_n)\in \Delta^d(S^1,n)$. Then there is a permutation $\sigma\in\sn$ bringing this configuration to a counterclockwise ordering $(x_{\sigma (1)},
\ldots, x_{\sigma (n)})$. Let $s_i$ be the arc distance (divided by $2\pi$) measured counterclockwise between $x_{\sigma (i)}$ and $x_{\sigma (i+1)}$. When $x_i\neq x_j$ for $i\neq j$, the choice of $\sigma$ is unique up to cyclic permutation and there is a well-defined map
\begin{eqnarray*}
F(S^1,n)&\lrar& (S^1\times \mathring{\Delta}_{n-1})\times_{\bbz_n}\sn\\
(x_1,\ldots, x_n)&\longmapsto&[(x_{\sigma (1)}, (s_1,\ldots, s_n)); \sigma ]
\end{eqnarray*}
which is a homeomorphism. Here $(s_1,\ldots, s_n)$ is in the open simplex
$\mathring{\Delta}_{n-1}$ if and only if none of the $s_i$'s is zero.
When there is collision, i.e.
$d>1$, then the choice of $\sigma$ up to cyclic permutation is not unique anymore but there
is a map at the level of unordered configuration spaces
\begin{eqnarray*}
B^d(S^1,n)&\lrar& S^1\ltimes_{\bbz_n} \Delta_{n-1}(d)\\
\ [x_1,\ldots, x_n]&\longmapsto&[x_{\sigma (1)}; (s_1,\ldots, s_n)]
\end{eqnarray*}
where $\sigma$ again is any permutation bringing $(x_1,\ldots, x_n)$ into cyclic ordering.
This  map is independent of the choice of $\sigma$ and it is a homeomorphism with inverse
$$[x_{\sigma (1)}; (s_1,\ldots, s_n)]\longmapsto [x_{\sigma (1)}, x_{\sigma (1)}e^{i2\pi s_1},
x_{\sigma (1)}e^{i2\pi (s_1+s_2)},\ldots, x_{\sigma (1)}e^{i2\pi (s_1+\cdots + s_{n-1})}]\ .$$
Note that when $x_i=x_{i+1}$ in the cyclic ordering, $s_i=0$ so the faces of $\Delta_{n-1}$ where the $s_i$'s vanish (consecutively) correspond to when points come together.
\end{proof}

\bpr\label{bds1} Identify $S^1=[0,1]/_\sim$. Then addition
$m: B^d(S^1,n)\lrar S^1$,
$$m ([x_1,\ldots, x_n])=x_1+ x_2+\cdots +x_n$$
is a bundle map with fiber $\Delta_{n-1}(d)$. In particular $m$ is a homotopy equivalence.
\epr

\begin{proof}
The composite
$$\rho : \xymatrix{S^1\times_{\bbz_n}\Delta_{n-1}(d)\ar[r]&B^d(S^1,n)\ar[r]^m&S^1
}$$
sends $(x, (s_1,\ldots, s_n))$ to $nx+ (n-1)s_1+(n-2)s_2 + \cdots + s_{n-1}$. This map
is well-defined on orbits since
$\rho (x+s_1, (s_2,\ldots, s_n,s_1)) = \rho (x, (s_1,\ldots, s_n))$.
The preimage of a point $y\in S^1$ under $m$ are all unordered tuples $[x_1,\cdots , x_n]$ such that $x_1+x_2+\cdots + x_n = y \mod \bbz$. All preimages are homeomorphic and we can choose $y=0$. The preimage $\rho^{-1}(0)$ consists of all classes $[x,(s_1,\ldots, s_n)]$ such that
$$(n-1)s_1+(n-2)s_2 + \cdots + s_{n-1} + nx \mod\bbz$$
We wish to show this is a copy of $\Delta_{n-1}(d)$. Consider the map $\phi : \Delta_{n-1}(0)\lrar \rho^{-1}(0)$ defined as follows.
Given $(s_1,\ldots, s_n)$, $\sum s_i=1$, let $m_s$ be the sum
${-1\over n}\left((n-1)s_1+(n-2)s_2 + \cdots + s_{n-1})\right)$ brought modulo $\bbz$ to the interval $[0,1]$ and define
$$\phi : (s_1,\ldots, s_n)\longmapsto [m_s, (s_1,\ldots, s_n   )]\in
S^1\times_{\bbz_n}\Delta_{n-1}(d)
$$
This map is well-defined and continuous. It is surjective by construction. It is also injective for the following reason. If
$s=(s_1,\ldots, s_n)$ and $s'=(s'_1,\ldots, s'_n)$ map to the same point under $\phi$, they must be
the same up to cyclic permutation. Let's assume $s' = (s_{k+1}, \ldots , s_n,s_1,s_2,\cdots, s_{k})$, $0<k< n$ ($s_0=s_n$).
A quick computation shows that
$$m_{s'} = m_s + s_1+\cdots +s_k- {k\over n}$$
But in  $S^1\times_{\bbz_n}\Delta_{n-1}(d)$,
$[m_{s'}, (s_1',\ldots, s_n')] = [m_s-{k\over n}, (s_1,\ldots, s_n)]$
so that $\phi (s')$ can never be $\phi (s)$ unless $k=0$ or $s_i=s'_i={1\over n}$. In both cases $s=s'$. This proves the injectivity and hence that $\phi$ is a homeomorphism. It remains to check that $\rho$ is a bundle map and this is left as an exercise.
\end{proof}

\bre (Morton) When $d=1$, $m: \sp{n}(S^1)\lrar S^1$ is an $n-1$-disk bundle which is trivial if and only if $n$ is odd. The open disk bundle is $B(S^1,n)$ and its sphere bundle is $B_2(S^1,n)$.
\ere

\subsection{Wedges of Circles} $B^d(S^1,n)\simeq S^1$ as discussed. The situation gets more complicated quickly for other graphs. The following is a neat little application of our constructions for the case $d=1$.

\bpr\label{wedgecircle} $\displaystyle B(\bigvee^kS^1,2)$ is homotopy equivalent to
$\displaystyle\bigvee^{{3\over 2}k(k-1)+ 1}S^1$.
\epr

\begin{proof}
Let's first understand the case $k=2$.

We will write $B(S^1\vee S^1,2)$ as the union of three subspaces:
$X_1=\{[(x,\ast ), (y,\ast)]\ |\ x\neq y\}$, then
$X_2=\{[(\ast, x),(\ast, y)]\ |\ x\neq y\}$ and
$X_3 = \{[(x,\ast), (\ast, y)]\ |\ (x,y)\neq (\ast,\ast)\}$. We have that
$$X_1\cong B(S^1,2)\ ,\  X_2\cong B(S^1,2)\ \ ,\ X_3\cong (S^1\times S^1)^* ,$$
where $(S^1\times S^1)^*$ means the punctured torus $S^1\times S^1-\{(\ast,\ast)\}$. Notice that $X_1\cap X_2=\emptyset$ while $X_1\cap X_3=\{(x,\ast ), (\ast,\ast), x\neq\ast\}\cong (S^1)^*$ are punctured circles hence contractible intervals. The punctured torus $X_3$ deformation retracts onto a wedge $S^1\vee S^1$. During this deformation both punctured circles corresponding to the intersection with $X_1$ and $X_2$ retract onto the wedgepoint. After the retraction we obtain a wedge $S^1\vee S^1\vee Y_1\vee Y_2$ where each $Y_i = X_i/_\sim$ is the open m$\ddot{\hbox{o}}$bius band $X_i=S^1\times ]0,1[$ with an interval $*\times ]0,1[$ retracted to a point. Therefore $Y_i\simeq S^1$ and the claim follows in this case.

\begin{figure}[htb]
\begin{center}
\epsfig{file=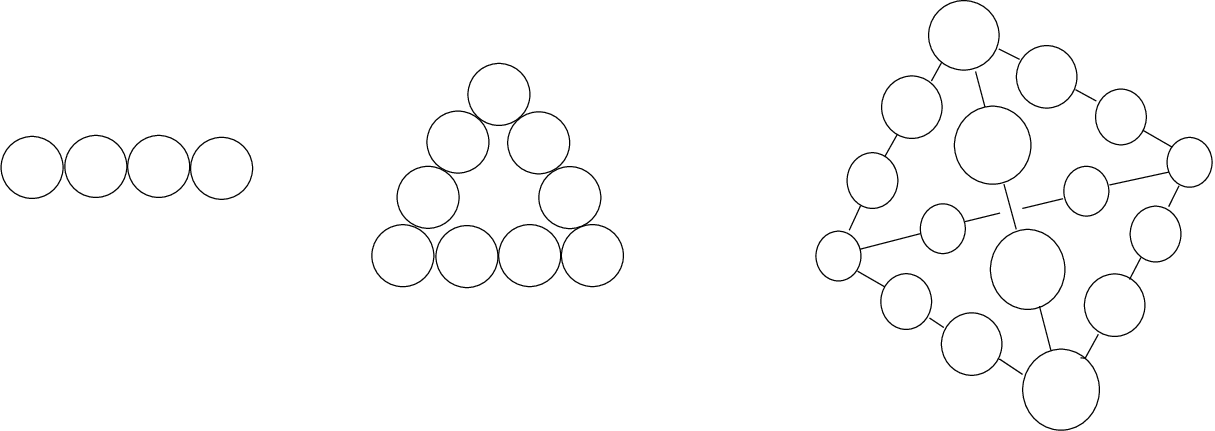,height=1.4in,width=3.8in,angle=0.0}
\caption{An intermediate homotopy type of $B(\bigvee^kS^1,2)$ for
$k=2,3$ and $4$ respectively. These are strings of circles making up a necklace in the shape of a $k-1$-dimensional simplex}\label{necklace}
\end{center}
\end{figure}

For the general case of a bouquet of $k$-circles, $k>2$, we write an element from the $i$-th leaf as $x^i$. Then $B(\bigvee^kS^1,2)$ becomes the union of subspaces
$X_{i,j}:=\{[(x^i,\ast ), (y^j,\ast)]\ |\ x^i\neq y^j\ \hbox{if}\ i=j\}$,
$X_i^j := \{[(x^i,\ast ), (\ast,y^j)]\ |\ (x^i, y^j)\neq (\ast,\ast)\ \hbox{if}\ i=j\}$ and
$X^{i,j}:= \{[(\ast, x^i ), (\ast,y^j)]\ |\ x^i\neq y^j\ \hbox{if}\ i=j\}$ over all $k\geq i\geq j\geq 1$.
As before $X_{i,i}= B(S^1,2)$ is the open m$\ddot{\hbox{o}}$bius band. For $i>j$,
$X_{i,i}$ and $X_{j,j}$ are disjoint. Also and as is clear,
$X_i^j\cap X_r^s = \emptyset$ if $\{i,j\}\neq \{r,s\}$. Each union
$B_{i,j}:=X_{i,j}\cup X_{i,i}\cup X_i^j$ is the sub-configuration space of $2$ points on the $i$-th and $j$-th leaves and hence is up to homotopy a wedge of $4$ circles. The homotopy deforming each $X_{i,i}$ to $S^1$ is the same if performed in $B_{i,j}$ or $B_{i,k}$. This is to say that the homotopies deforming $B_{i,j}$'s to a wedge of $4$ circles are compatible and we obtain a deformation retract of $B(\bigvee^kS^1,2)$ which looks like a \textit{necklace} of circles tied in the shape of the $k-1$-dimensional simplex. This is depicted in figure \ref{necklace} for $k=2,3$ and $4$. The homotopy type of this space is not hard to work out: it is a wedge of all those circles appearing in the necklace with another wedge of circles describing the homotopy type of the $1$-skeleton of $\Delta_{k-1}$. In the necklace there is one circle for each vertex of the $k-1$-simplex and two circles for each edge, this gives a total of $k^2$ circles. On the other hand the one-skeleton of the $k-1$-simplex, denoted by $\Delta_{k-1}^{(1)}$, is homotopy equivalent to $\bigvee^{N}S^1$ where $N = {1\over 2}k(k-3) + 1$ circles. Indeed
the Euler characteristic
$$\chi (\Delta_{k-1}^{(1)}) = \#\hbox{edges}-\#\hbox{vertices} =
{1\over 2}k(k-1) - k = {1\over 2}k(k-3) ,
$$
and this must be $\chi (\bigvee^NS^1) = N-1$.
Putting this together yields
$$B(\bigvee^kS^1,2)\simeq \bigvee^{k^2}S^1\vee\bigvee^{{1\over 2}k(k-3) + 1}S^1\simeq
\bigvee^{{3\over 2}k(k-1)+ 1}S^1 $$
and the proof is complete.
\end{proof}

\bre The first homology group of $B(\Gamma, n)$ for graphs has been worked out in \cite{kp}. Their method uses discrete Morse theory. In particular one can deduce from their Theorem 3.16 that $H_1(B(\bigvee^kS^1,2))=\bbz^{1+3k(k-1)/2 }$ in full agreement with our Proposition \ref{wedgecircle} [in their theorem one uses that the braid index is $2$, $N_1 = {2k(k-1)-{k(k-1)\over 2}-(k-1)}$ and the first betti number of the graph is of course $k$]. In the case of trees $T$, the homology groups of the unordered configuration space $B(T,n)$ are torsion free and their ranks computed by Farley (see \cite{kp} as well).
\ere


\section{The Localization Principle and the Case of Manifolds}\label{main1b}

Our main approach is to find conditions on $X$ so that the inclusion $B^d(X,n)\hookrightarrow \sp{n}(X)$ induces an isomorphism on some homotopy groups through a range.

 We start with a preliminary lemma. We say a space $X$ is \textit{locally punctured connected} if for every $x\in X$ and  neighborhood $U$ of $x$, there is an open $V$, $x\in V\subset U$ such that $V-\{x\}$ is connected.

\ble\label{connected} Let $X$ be path-connected which is not a point, locally contractible. If $d\geq 2$, then both $\Delta^d(X,n)$ and $B^d(X,n)$ are connected. If furthermore $X$ is locally punctured connected, then both $\Delta^d(X,n)$ and $B^d(X,n)$ are connected for all $d\geq 1$.
\ele

\begin{proof}
For both claims, it suffices to show that $\Delta^d(X,n)$ is connected. We need join $(x_1,\ldots, x_n)$ to $(y_1,\ldots, y_n)$ by a path, for any two choices of tuples in $\Delta^d(X,n)$. By deforming locally, we can arrange that the $x_i$'s and the $y_j$'s are all pairwise distinct. Now $X$ is path-connected so there is a path $\gamma_i$ from $x_i$ to $y_i$. Via $\gamma_1$ we construct a path in $\Delta^d(X,n)$ from
$(x_1,x_2,\ldots, x_n)$ to $(y_1,x_2,\ldots, x_n)$ by putting $\gamma_1(t)$ in the first coordinate. At any given time $t\in [0,1]$, $\gamma_1(t)$ can only coincide with one $x_i$ at a time, and hence this path is well-defined in $\Delta^d(X,n)$ if $d\geq 2$. Construct next the path from
$(y_1,x_2,\ldots, x_n)$ to $(y_1,y_2,x_3,\ldots, x_n)$ by putting
$\gamma_2(t)$ in the second coordinate. This is again a well defined path in $\Delta^d(X,n)$. We can continue this process. The composition $\gamma_n\circ\cdots\circ\gamma_1$ is a path in $\Delta^d(X,n)$ from
$(x_1,\ldots, x_n)$ to $(y_1,\ldots, y_n)$.

To establish the second claim, we proceed by induction on $n\geq d$. For $n=d$, $\Delta^d(X,d)=X^d$ and there is nothing to prove. For $n>d$, consider the projection
$$\Delta^d(X,n)\lrar \Delta^d(X,n-1)$$
which omits the last coordinate. The preimage of a tuple $(x_1,\ldots, x_{n-1})$ is $X-\{x_{i_1},\ldots, x_{i_j}\}$ if $x_{i_r}$ repeats $d$-times
in the tuple. Since $X$ is locally punctured connected, this preimage is connected. Since the base space of this projection is also connected by inductive hypothesis, it follows that the total space is connected as desired.
\end{proof}

For the higher homotopy groups, the starting point is the following principle which relates the local connectivity properties of a space to its global properties. All spaces appearing below are connected. The following result is in \cite{kloeckner}, Theorem 1.4.

\bth\label{localizing}(Localization)
Let $X$ be a Hausdorff topological space and $Y$ be a closed subset of $X$. If for every point $y\in Y$, and every neighborhood $U\subset X$ of $y$, there is an open $V\subset U$ containing $y$ such that the pair $(V, V\setminus  Y)$ is $k$-connected, $k\geq 0$, then the pair $(X, X\setminus Y)$ is $k$-connected.
\end{theorem}

We recall what it means for the pair $(X,A)$ to be $k$-connected or that $\pi_i(X,A) =0$ for all $i\leq k$ (\cite{hatcher}, chapter 4). If $k\geq 1$, this means that every map $(I^r,\partial I^r) \lrar (X,A)$ from the closed cube $I^r$, $1\leq r\leq k$, is homotopic (relative its boundary) to a map $I^r\lrar A$. For any $x\in X$, $(X,\{x\})$ is $1$-connected if and only if $X$ is simply connected. Being $0$-connected, or equivalently writing
$\pi_0(X,A)=0$, means in our terms that $X$ and $A$ are connected and that any point in $X$ is connected by a path to a point in $A$. Note that in the Theorem, if either $V\backslash Y$ or $X\backslash Y$ is not connected, then the Theorem fails.

\bex $X=\bbr^3$ and $L$ a line in $\bbr^3$. The pair $(X,X\setminus L)$ is $1$-connected but not $2$-connected. Indeed take a square which is intersected transversally through its interior by $L$. That square cannot be deformed away from $L$ with the boundary being kept fixed.
\eex

The following is a consequence of Theorem \ref{localizing}. We say that a closed subset $Y$ in $X$ is \textit{tame} if there is a neighborhood $N$ of $Y$ so that $N$ deformation retracts onto $Y$ and $X\setminus Y$ deformation retracts onto $X\setminus N$. Submanifolds are tame and so are subcomplexes of simplicial complexes (Proposition \ref{munkres}).

\bco\label{localglobal} Let $Y$ be a tame subspace of $X$ and suppose for every
$y\in Y$ and neighborhood $U$ of $y$ in $X$, there is a contractible neighborhood $V\subset U$, such that $Y\cap V$ is tame in $V$ and $V\setminus Y$ is $k$-connected, $k\geq 0$. Then $\pi_i(X)\cong\pi_i(X\setminus Y)$
for $i\leq k$.
\eco

\begin{proof} The point is that when $Y$ is tame in $X$,  Theorem \ref{localizing} implies that the induced map $\pi_{k}(X\backslash Y) \rightarrow \pi_{k}(X)$ is surjective, and $\pi_i(X\backslash Y)\rightarrow\pi_i(X)$ is an isomorphism for $i\leq k-1$. Let's show that for $(V,y)$ as in the statement of the theorem, the pair $(V, V\setminus Y)$ is $k+1$-connected. Since $V\cap Y$ is tame in $V$, choose a neighborhood $N$ of $Y$ in $V$ that deformation retracts onto $Y$ and such that $V\setminus N$ deformation retracts onto $V\setminus Y$. We can replace up to homotopy the pair
$(V,V\setminus Y)$ by $(V,V\setminus N)$ where now $V\setminus N$ is closed in $V$.
We can apply the long exact sequence in homotopy of the pair $(V, V\setminus N)$
$$ \rightarrow \pi_{k+1}V\longrightarrow  \pi_{k+1}(V,V\setminus N)\fract{\partial}{\longrightarrow } \pi_{k}(V\setminus N)\longrightarrow \cdots\rightarrow  \pi_1(V,V\setminus N)\fract{\partial}{\longrightarrow } \pi_0(V\setminus N)\rightarrow \pi_0(V)$$
Since for $i\leq k$, $\pi_i(V\setminus N)=0=\pi_{i+1}(V )$, we see that $
\pi_i(V,V\setminus Y)\cong \pi_i(V,V\setminus N) = 0$ for $i\leq k+1$. From Theorem \ref{localizing} it follows that $(X,X\setminus Y)$ is $k+1$-connected. The same argument as above with the long exact sequence of the pair $(X,Y)$ with $Y$ tame in $X$ shows that $\pi_i(X)\cong\pi_i(X\setminus Y)$
for $i\leq k$.
\end{proof}

\bre\label{localglobal2} In the case of a submanifold $S$ in $M$ of codimension $d$, a neighborhood of a point deformation retracts onto a sphere $S^{d-1}$ which is $d-2$-connected. By the previous corollary this gives that $M$ is weakly equivalent to $M-S$ up to dimension $d-2$ (Proposition \ref{manifoldcase}). A similar argument applies when $S=\bigcup S_j$ is the union of submanifolds intersecting transversally.
\ere

The following key Lemma shows how we can apply the above results to diagonal arrangements.

\ble\label{key} Let $X$ be a finite simplicial complex so that for every $x\in X$ and neighborhood $U$ of $x$, there is a subneighborhood $V$ containing $x$ such that $\Delta^d(V,k)$ (respectively $B^d(V,k)$) is $r$-connected for any $k\geq 1$. Then $\pi_i(\Delta^d(X,n))\cong\pi_i(X^n)$ (resp. $\pi_i(\sp{n}X)\cong\pi_i(B^d(X,n))$ for $i\leq r$.
\ele

\begin{proof}
We have to estimate the connectivity of the pair $(X^n,\Delta^d(X,n)) = (X^n, X^n-\Delta_{d+1}(X,n))$ (resp. that of $(\sp{n}M,\sp{n}(M)-B_{d+1}(M,n))$. Note that $\Delta_{d+1}(M,n)$ (resp. $B_{d+1}(M,n)$) is tame in $M^n$ (resp. $\sp{n}M$) according to \S\ref{compactification}. One can check they verify the hypothesis of Corollary \ref{localglobal}. In the ordered case, choose a point in $\Delta_{d+1}(X,n)$ which after permutation can be brought to the form
\begin{equation}\label{thepoint}
(\underbrace{x_1,\cdots, x_1}_{i_1}, \underbrace{x_2,\cdots,x_2}_{i_2},
\ldots, \underbrace{x_r,\cdots ,x_r}_{i_r})\ .
\end{equation}
with $x_i\neq x_j$ if $i\neq j$, $\sum i_\alpha = n$ and $i_1>d$.
A neighborhood $W$ of this point in $X^n$ is homeomorphic to
$V_1^{i_1}\times\cdots\times V_r^{i_r}$ where $V_i$ is a contractible neighborhood of $x_i$ in $X$, the $V_i$'s pairwise disjoint. Clearly
$$W-\Delta_{d+1}(X,n)\cong \Delta^d(V_1,i_1)\times\cdots\times \Delta^d(V_r,i_r)\ .$$
By hypothesis we can assume all the $\Delta^d(V_i,i_j)$ to be $r$-connected
so that $W-\Delta_{d+1}(X,n)$ is also $r$-connected and hence, by Corollary \ref{localglobal}, $\pi_i(\Delta^d(X,n))=\pi_i(X^n-\Delta_{d+1}(X,n))\cong\pi_i(X^n)$ for $i\leq r$.

A similar proof holds in the unordered case. Given a point in $B_{d+1}(M,n)\subset\sp{n}(M)$ as in (\ref{thepoint}), a small contractible neighborhood of it in $\sp{n}M$ is
$U \cong \sp{i_1}(V_1)\times\sp{i_2}(V_2)\times\cdots\times\sp{i_r}(V_r)$, the $V_i$'s pairwise distinct, and
\begin{equation}\label{theneighborhood}
B^d(U,n)= U-B_{d+1}(X,n)\cong
B^d(V_1,i_1)\times
B^{d}(V_2,i_2)\times\cdots\times B^d(V_r, i_r)\ .
\end{equation}
If we choose each $V_j$ so that $B^d(V_j,i_j)$ is $r$-connected (hypothesis), the complement $B^d(U,n)$ will also be $r$-connected and the claim follows again from Corollary \ref{localglobal}.
\end{proof}

In the case of a manifold we can already make the following easy conclusions.

\bco\label{manifold1}
Let $M$ be a manifold of dimension $m\geq 1$. \\
(i) If $m d\geq 3$ then
$\pi_i(\Delta^d(M,n))\cong \pi_i(M)^n$ for $i\leq md-2$. \\
(ii) If $m\geq 2$ and $d\geq 2$, then
$\pi_1(B^d(M,n))\cong H_1(M,\bbz)$.
\eco

\begin{proof}
Every point of $M$ has a neighborhood homeomorphic to $\bbr^m$.
The fat diagonal $\Delta_{d+1}(\bbr^m,n)$ in $(\bbr^m)^n$ has codimension $mn-m(n-d) = md\geq 3$ so its complement $\Delta^d(\bbr^m,n)$ is $md-2$-connected (Proposition \ref{manifoldcase}). Now apply Lemma \ref{key} to get (i). On the other hand  $\Delta^d(\bbr^m,k)$ is simply connected if $d\geq 2$ and $m\geq 2$, so by Armstrong's result (Corollary \ref{armstrong}), $\pi_1(B^d(\bbr^m,k))$ is also trivial and (ii) follows.
\end{proof}

\bre As we pointed out, Corollary \ref{manifold1} (i) is not true for $md=2$ as illustrated by $F(\bbr^2,2)\simeq S^1$. This corollary is a special case of Theorem \ref{main1}.
Also let's point out that $\Delta^d(\bbr^m,n)$ has torsion free homology starting with spherical
classes in $dm-1$ as already indicated and all homology classes being represented by products of spheres
\cite{dt}.
\ere

We now derive Theorem \ref{main1} when $X$ is a manifold.
Again $X$ is $r$-connected if $\pi_i(X)=0$ for $0\leq i\leq r$.

\ble\label{connectivity} Let $\Omega^m_*(-)$ denote a connected component of the loop space $\Omega^m(-)$, $m\geq 1$ and $d\geq 1$. Then $\Omega^m_*\sp{d}S^m$ is $2d-2$ connected.
\ele

\begin{proof} Let's review the simplest cases.
The case $d=1$ is obvious since $\Omega^mS^m$ breaks down into components indexed by the integers, and each component is $0$-connected but not $1$-connected since $\pi_1(\Omega^m_*S^m)\cong \pi_{m+1}(S^m)$ is $\bbz$ if $m=2$ and $\bbz_2$ if $m\geq 3$. When $m=1$,
$\sp{d}S^1\simeq S^1$ so that $\Omega_*S^1$ is contractible and hence certainly $2d-2$ connected for any $d$. When $m=2$, $\sp{d}S^2\cong{\mathbb P}^d$ is complex projective space and
$$\Omega^2\sp{d}S^2 = \Omega^2{\mathbb P}^d\cong {\mathbb Z}\times\Omega^2S^{2d+1} .$$
Each component is a copy of $\Omega^2 S^{2d+1}$ which is $2d-2$-connected and the bound is sharp.

In general we invoke Theorem 5.9 of \cite{kt} which states that for $r$-connected $X$, $r\geq 1$,
\begin{equation}\label{connectivity2}
\pi_i(\sp{n}X)\cong\tilde H_i(X;\bbz )\ \ \ ,\ \ \ 0\leq i\leq r+2n-1\ .
\end{equation}
This gives that for $i\geq 1$ and $m\geq 2$,
$$\pi_i(\Omega^m_*\sp{d}S^m)\cong\pi_{i+m}(\sp{d}S^m)\cong
H_{i+m}(S^m) = 0\ \ \ ,\ \ \ i+m\leq (m-1)+2d-1 = m+2d-2\ .$$
This gives $i\leq 2d-2$ and a lower bound for the connectivity is $2d-2$.
\end{proof}

\bpr\label{bdrmn} Assume $m\geq 2, n\geq d\geq 1$. Then $B^d(\bbr^m,n)$ is $2d-2$-connected. Moreover
if $X$ is a $1$-connected manifold and $n\geq 2$, then $\pi_i(B^d(X,n))\cong \tilde H_i(X;\bbz )$ for $0\leq i\leq 2d-2$.
\epr

\begin{proof} This relies on results from \cite{kallel,kt}. The case $d=1$ being trivial,
we assume $d\geq 2$. Consider the sequence of embeddings
\begin{equation}\label{themap}
\tau_n: B^d(\bbr^m,n)\hookrightarrow B^d(\bbr^m,n+1)\ \ \,\ \ \
[x_1,\ldots,x_{n}]\mapsto [x_1,\ldots, x_{n},|x_1|+\cdots +|x_n|+1]\ .
\end{equation}
The direct limit is $B^d(\bbr^m,\infty)$ and it is shown in \cite{kallel} that there is a ``scanning map''
$$\tau : B^d(\bbr^m,\infty)\lrar\Omega^m_*\sp{d}S^m$$
which induces a homology isomorphism. Since both spaces are simply connected when $d\geq 2$ (Lemma \ref{armstrong} and Lemma \ref{connectivity}) of the homotopy type of CW complexes, the map $\tau$ is a homotopy equivalence.
Moreover, the maps $\tau_n$ in (\ref{themap}) induce homology embeddings according to
(\cite{zanos}, chapitre 3). Iterating we get homology embeddings
$$H_*(B^d(\bbr^m,d+1))\hookrightarrow H_*(B^d(\bbr^m,n)\hookrightarrow
H_*(B^d(\bbr^m,\infty))\cong H_*(\Omega^m_*\sp{d}S^m)\ .$$
By Lemma \ref{connectivity} the groups on the extreme right are trivial for
 $*\leq 2d-2$.  This gives that $H_*(B^d(\bbr^m,n)) = 0$ for $n\geq d$ and $*\leq 2d-2$. Since the space is simply-connected, it is $2d-2$-connected as well. It then follows by Lemma \ref{key} that $\pi_i(B^d(X,n))\cong\pi_i(\sp{n}(X))$ for $i\leq 2d-2$. This proves the main statement. In the case $X$ is $r$-connected with $r\geq 1$, it follows by the inequality in (\ref{connectivity2}), since $2d-2\leq r+2n-1$, that
$\pi_i(B^d(X,n))\cong\pi_i(\sp{n}(X))\cong \tilde H_i(X;\bbz )$ in the range
of dimensions $0\leq i\leq 2d-2$.
\end{proof}

\bex Consider the case $B^{n-1}(S^2,n)$, $n\geq 3$. Since $\sp{n}(S^2)\cong\bbp^n$ is a $2n$-dimensional manifold, by Proposition \ref{manifoldcase}, $\pi_i(B^{n-1}(S^2,n))\cong \pi_i(\bbp^n)$ for
$1\leq i\leq 2(n-1)-2=2n-4$. On the other hand from the Hopf fibration,
$\pi_i(\bbp^n)\cong\pi_i(S^{2n+1})$ for $i>2$ and $\pi_2(\bbp^n)=\bbz$. This shows precisely that $\pi_i(B^{n-1}(S^2,n))\cong H_i(S^2,\bbz)$ for $1\leq i\leq 2n-4$ as expected.
\eex

The claim that $B^d(\bbr^k,n)$ is $2d-2$-connected has an alternative nice proof in the case $d=n-1$.

\ble\label{bn-1rk} $B^{n-1}(\bbr^k,n)$ is $2n-4$-connected, $n\geq 2, k\geq 1$.
\ele

\begin{proof} The case $k=1$ is trivial. We let $k\geq 2$ and invoke some main results from \cite{kk,kt}. Let $S$ be the unit sphere as in Lemma \ref{trivial} and let $Q_{n,k}$ be its quotient under the $\sn$-action. We have already indicated that $Q_{n,k}\simeq B^{n-1}(\bbr^k,n)$.
On the other hand, according to Theorems 1.1, 1.3 and 1.5 of \cite{kk}, we have that
\begin{equation}\label{identification1}
\Sigma^{k+1}Q_{n,k}\simeq \bsp{n}(S^{k}) ,
\end{equation}
where $\Sigma$ means suspension and $\bsp{n}(Y)$ means the ``symmetric smash'' $Y^{\wedge (n)}/\sn$, which is also the cofiber of the embedding of
$\sp{n-1}Y$ into $\sp{n}Y$ induced by adjoining a basepoint to an unordered tuple $[x_1,\ldots, x_{n-1}]$. It is shown (\cite{kk}, Theorems 1.2 and 1.3) that if $X$ is $r$-connected, then $\bsp{n}(\Sigma X)$ is $2n+r-1$-connected. This gives that $\bsp{n}(S^{k})=\bsp{n}(\Sigma S^{k-1})$ is $2n+k-3$ connected, and hence so is $\Sigma^{k+1}Q_{n,k}$ by (\ref{identification1}). Since in this range
$Q_{n,k}$ is already simply connected, it must therefore be $2n-4$-connected.
\end{proof}

\bre That the connectivity bound in the above theorem doesn't depend on $k$ is not surprising. Indeed when $n=2$, $B(\bbr^k,2)\simeq\bbr P^{k-1}$ and this is never $1$-connected no matter what $k$ is.\ere

\section{An Equivariant Deformation Retract of Diagonal Complements}\label{compactification}

Let $X_\ast$ be an abstract simplicial complex and $|X_\ast|$ its geometric realization. Let $A_\ast$ be a subcomplex of $X_\ast$.
We say a subcomplex $A_\ast$ of $X_\ast$ is \textit{full} if every simplex of $X_*$ whose vertices are in $A_\ast$ is itself in $A_\ast$. The following fundamental result (called the "retraction lemma" in \cite{blago2}) can be found in Munkres' book \cite{munkres}, Lemma 70.1.

\bpr\label{munkres} Let $A_\ast$ be a full subcomplex of the finite simplicial complex $X_\ast$. Let $C_\ast$ consist of all simplices of $X_\ast$ that are disjoint from $A_*$. Then $|A_\ast|$ is a deformation retract of $|X_\ast|-|C_\ast|$, and $|C_\ast|$ is a deformation retract of $|X_\ast|-|A_\ast|$.
\epr

The argument of proof is short but instrumental to extract useful properties of this ``compactification". We review this argument.
The fact that $A_*$ is full says that $C_*$ is also full, and that simplexes of $X_*$ consist of simplexes in $C_*$, simplexes in $A_*$ and simplexes of the form
$$\sigma *\tau  \ \ \ ,\ \ \sigma\in A_*\ ,\ \tau\in C_* ,$$
where $\sigma *\tau$ is the join of both simplexes. The following figure illustrates the situation when $X_*$ is the full simplex $\Delta_3$ on $4$ vertices $v_0,v_1,v_2,v_3$, $A_*=[v_0v_1]$ and $C_*=[v_2v_3]$

\begin{figure}[htb]
\begin{center}
\epsfig{file=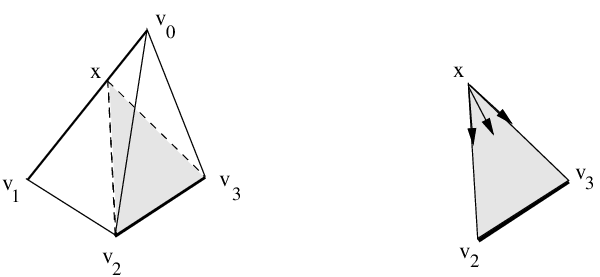,height=1.5in,width=3.5in,angle=0.0}
\caption{Munkres' deformation along the join (right) after deleting the apex $x$}\label{deformation}
\end{center}
\end{figure}

The deformation of $|X_*|-|A_*|$ onto $|C_*|$ is as in the figure. It
starts at a point $tx + \sum_{i\in I} s_iv_i$, with $v_i$ vertices in $C_*$, $i\in I$, $t+\sum s_i=1$, $t\neq 1$, and ends at the point
$\sum {t_j\over\sum s_i}v_j, j\in I$.

Two important consequences are in order:

\begin{itemize}
\item If $A_\ast$ is full, $|X_\ast|-|A_\ast|$ deformation retracts onto the \textit{largest} subcomplex that does not meet $|A_\ast|$. Note that if $A_\ast$ is not full, then its first barycentric subdivision $SdA_\ast$ is always full in $SdX_\ast$. The barycentric subdivision comes with a natural ordering on vertices.
\item The deformation retraction illustrated in figure \ref{deformation} has the property that if it starts in a simplex of $X_\ast$ it will stay in that simplex (and deforms onto a face of it).
\end{itemize}

For ease we will write $X$ for either $X_\ast$ or its realization. The context will be clear.

Munkres' observation nicely applies to the diagonal arrangements. Given $X$ an ordered simplicial complex, then $X^n$ can be given naturally a structure of a simplicial complex such that the various diagonals are subcomplexes (see \cite{nakaoka}, \S1, also proof of Lemma \ref{property1} below).
We can then apply Proposition \ref{munkres} to the configuration space
$X^n-\Delta_{d+1}(X,n)$. Among all diagonal arrangements, only the thin diagonal $\Delta_{n}(X,n)$ is full. We therefore have to pass to a barycentric subdivision.
Let $Sd(X^n)$ be the barycentric subdivision of $X^n_*$. This restricts to
$Sd(\Delta_{d+1}(X,n))$.

\ble\label{property1} There is $\sn$-equivariant deformation retraction of $\Delta^d(X,n)$ onto the largest subcomplex $W^d(X,n)$ not intersecting $|Sd(X^n)|-|Sd(\Delta_{d+1}(X,n))|$.
\ele

\begin{proof} That the complement deformation retracts onto
$W^d(X,n)$ is a direct consequence of Proposition \ref{munkres} as applied to the pair
$(Sd(X^n),Sd(\Delta_{d+1}(X,n)))$ with $Sd(\Delta_{d+1}(X,n))$ being full. We need check this deformation is equivariant under the symmetric group action. Recall that the simplicial decomposition of $X^n$ is made out as follows, where $X$ of course is an \textit{ordered} simplicial complex \cite{nakaoka}. A vertex of $X^n$ is of the form $(v_1,\ldots, v_n)$ where $v_i$ is a vertex of $X$. Different $(q+1)$-vertices
$$w_0=(v_{01},\ldots, v_{0n})\ , \ w_1=(v_{11},\ldots, v_{1n})\ \cdots \
w_0=(v_{q1},\ldots, v_{qn})$$
form a $q$-dimensional simplex if and only if for each $k=1,2,\ldots, n$, $(q+1)$-vertices $v_{0k},v_{1k},\dots, v_{qk}$ are contained in a simplex of $X$ and it holds that
$v_{0k}\leq v_{1k}\leq\cdots\leq v_{qk}$ (see figure \ref{decomposition} for the decomposition of $X^3$ in the case $X=[0,1]$ with vertices $[0]\leq [1]$). Notice as asserted that $\Delta_2([0,1],3)$ is not full since the two-simplex (bottom) $([0,0,0],[1,0,0],[1,1,0])$ has all three vertices in $\Delta_2(X,3)$ but it is not itself a simplex of $\Delta_2(X,3)$.

\begin{figure}[htb]
\begin{center}
\epsfig{file=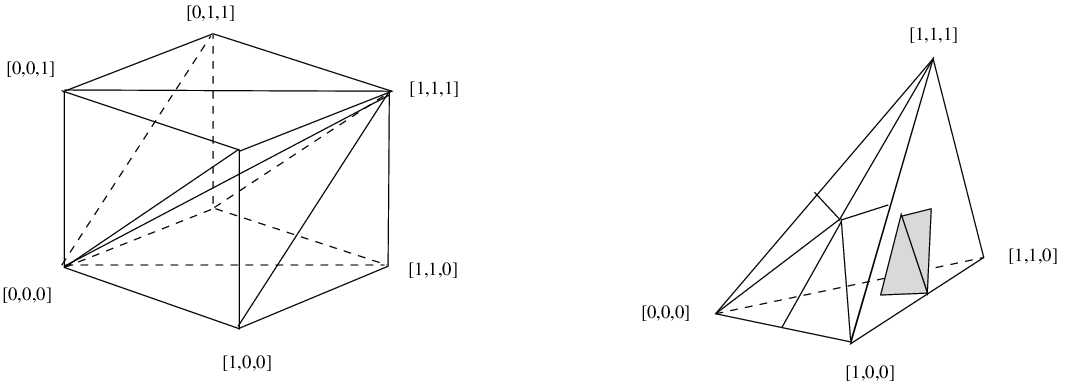,height=1.8in,width=4.5in,angle=0.0}
\caption{Left: simplicial decomposition for $[0,1]^3$ with $8$ vertices, $19$ edges, $18$ triangular faces and $6$ tetrahedral faces. Note $\Delta_2([0,1],3)$ is not full and we need pass to a barycentric subdivision. Right: the configuration space $|Sd(X^3)|-|Sd(\Delta_2(X,3)|$  deformation retracts onto the subcomplex $W^2([0,1],3)$ made out of
$6$ contractible connected components. The figure shows one such component in one tetrahedral face.}\label{decomposition}
\end{center}
\end{figure}

Generally a vertex is in $\Delta_{d+1}(X,n)$ if and only if it is of the form $(v_{1},\ldots, v_{n})$ for some vertices $v_1,\ldots, v_n$ of $X$ with
$v_{i_0}=\cdots = v_{i_{d}}$ for some choice of sequence $i_0<i_1<\cdots <i_d$.
Obviously every permutation acting on $X^n$ permutes vertices of $X^n_*$ and the order between them so it must take simplexes to simplexes. The action is simplicial and the quotient space $\sp{n}(X)$ inherits a cellular decomposition. Moreover the action remains simplicial after passing to a barycentric subdivision. Indeed since any new  introduced vertex is of the form ${1\over k}\sum v_i$, it is sent by $\sigma\in\sn$ to ${1\over k}\sum\sigma (v_i)$ which is the barycenter of $(\sigma (v_1),\ldots , \sigma (v_k))$.

After one subdivision, a simplicial neighborhood of $Sd(\Delta_{d+1}(X,n))$ consists of all simplices of $Sd(X^n)$ having at least one vertex of the form $(v_1,\ldots, v_n)$ with $v_{i_0}=\cdots = v_{i_{d}}$ for some sequence $i_0<i_1<\cdots < i_{d}$. This simplicial neighborhood is therefore $\sn$-invariant and its complement $W^d(X,n)$ is invariant as-well. Clearly the permutation action on $X^n$ commutes with Munkres` deformation since it takes combinations $\sum t_iv_i$ to $\sum t_i\sigma (v_i)$ (see Figure \ref{deformation}). It therefore descends to a deformation retraction of $B^d(X,n)$ onto $W^d(X,n)/\sn =: {\mathcal W}^d(X,n)$.
\end{proof}

\bco\label{mathcalw} For $X$ a finite simplicial complex, the $\sn$-quotient ${\mathcal W}^d(X,n)$ of $W^d(X,n)$ is a compact deformation retract of $B^d(X,n)$.
\eco

We need one more observation.

\ble\label{property2} Let $A$ be a subcomplex of $X$. The deformation retraction of $|Sd(X^n)|-|Sd(\Delta_{d+1}(X,n))|$ onto its compactified space
$W^d(X,n)$ restricts to a deformation retraction of
$|Sd(A^n)|- |Sd(\Delta_{d+1}(A,n))|$ onto $W^d(A,n)$.
\ele

\begin{proof} Since $A$ is a subcomplex of $X$,
$\Delta_{d+1}(A,n)$ is a subcomplex of $\Delta_{d+1}(X,n)$ and
$Sd(A^n)$ is a subcomplex of $Sd(X^n)$. Both $Sd(\Delta_{d+1}(X,n))$
and $Sd(\Delta_{d+1}(A,n))$
are full subcomplexes.
The assertion now follows from the fact that if the deformation retraction starts in a simplex of $Sd(X^n)$; in particular in $Sd(A^n)$, it will stay in that simplex.
\end{proof}


\section{Proof Theorem \ref{main1}}\label{proofofmain1}

We appeal to the following useful theorem of Steve Smale which is a generalization of classical results of Begle and Vietoris. A similar statement for maps between simplicial complexes can be deduced from work of Dror Farjoun  (\cite{dror}, Corollary 9.B.3, page 163).

\bth\label{smale}\cite{smale}
Let $X$ and $Y$ be connected, locally compact, separable metric spaces, and let $X$ be locally contractible. Let $f$ be a mapping of $X$ into $Y$ for which $f^{-1}$ carries compact sets into compact sets. If, for each $y\in Y$, $f^{-1}(y)$ is locally contractible and $r$-connected, then
the induced homomorphism $\pi_k(X)\lrar\pi_k(Y)$ is an isomorphism for $0\leq k\leq r$, and is onto  for $k=r+1$. \end{theorem}

Theorem \ref{smale} uses maps that are proper and preimages that are at least connected. Maps between configuration spaces obtained by projections are seldom proper. Combining the above theorem with \S\ref{compactification} yields however the following main result.

\bth\label{most}
Let $X$ be a connected finite simplicial complex with at least two vertices, $d\geq 2, n\geq 2$. Then
$$\pi_i(B^d(X,n))\cong \pi_i(\sp{n}(X))\ \ ,\ \ 0\leq i\leq 2d-2\ .$$
\end{theorem}

\begin{proof} The starting point is Lemma \ref{key} where it suffices to show that
$\pi_i(B^d(V,n))=0$ for $i\leq 2d-2$ for $V$ a small contractible neighborhood of a point in $X$. A neighborhood $V$ of $x\in X$ is of three types; either (i) Euclidean space, (ii)  half-space or (iii) it is a union of such half-spaces along a ``shared boundary".
 See figure \ref{contractnhbd}.
 \begin{figure}[htb]
\begin{center}
\epsfig{file=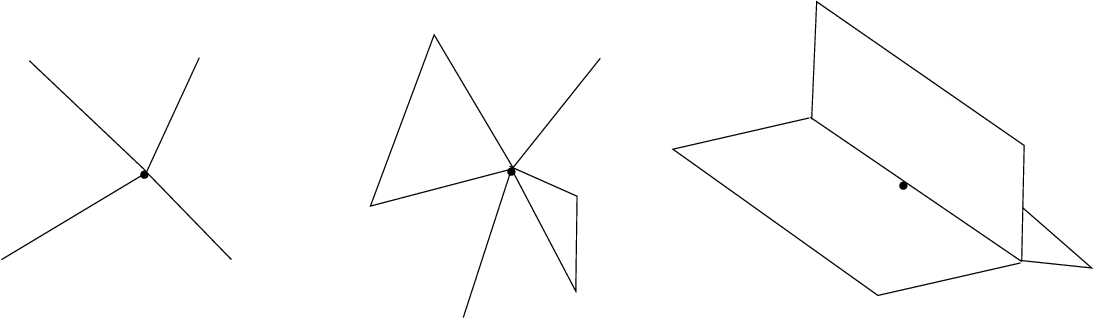,height=0.9in,width=3.8in,angle=0.0}
\caption{contractible neighborhoods of the dotted point in simplical $X$}\label{contractnhbd}
\end{center}
\end{figure}

We claim that in all cases, $B^d(V,n)$ is $2d-2$-connected.

In the case $x$ is an interior point of a simplex that is not a face of a larger simplex, it has a neighborhood $V\cong\bbr^m$ with $m\geq 1$. When $m=1$, $B^d(\bbr^1,n)$ is contractible. When $m\geq 2$,
$B^d(\bbr^m,n)$ is $2d-2$-connected according to Proposition \ref{bdrmn}.

If $x$ belongs to a boundary face, then $V$ is homeomorphic to half-space $H$ (with boundary). This half-space can be isotoped into its interior $\mathring{H}$ so we have a map  $B^d(H,n)\lrar B^d(\mathring{H},n)$ obtained from a deformation retraction (setting $t=1$). Since $B^d(\mathring{H},n)$ is $2d-2$-connected as seen prior, it follows immediately that $B^d(H,n)$ has the same connectivity (at least).

In the third and final case, $x$ lies in the intersection of two or more simplexes of $X$ as in figure \ref{contractnhbd}. Let $V$ be a contractible neighborhood made out of simplices which meet along a simplex $A$. Let $\Gamma$ be a simplex in $V$ of dimension $m$. Of course $A$ is in the boundary of $\Gamma$. Let
$$B^d(\Gamma, A,n) = \coprod_{0\leq k\leq n} B^d_A(\Gamma, k)/_{\sim}$$
where $B^0_A(\Gamma,k)=\ast$ is a given point in $A$
and $B^d_A(\Gamma, k)=B^d(\Gamma,k)\cup\sp{k}(A)$ (i.e the only points that can repeat more than $d$ times in $\Gamma$ are those that are in $A$).
The equivalence relation $\sim$ is such that
$x\sim\ast$ if $x\in A$ and $[x_1,\ldots,x_i,\ldots, x_k]\sim [x_1,\ldots, \hat x_i,\ldots, x_k]$ if $x_i\in A$. Here as customary $\hat x_i$ means the $i$-th entry is suppressed.
We have a projection
\begin{equation}\label{lambdamap}
\lambda : B^d (V,n)\lrar B^d (\Gamma, A,n),
\end{equation}
which sends a tuple $[x_1,\ldots ,x_n]$ to the new tuple obtained by replacing all $x_i\not\in\Gamma$ by $\ast$. On can view $\lambda$ as projection of $[x_1,\ldots ,x_n]$
 to the subtuple made up of those entries $x_i\in \Gamma$. This map is continuous by the very nature of the construction $B^d (\Gamma, A,n)$ (i.e. any entry $x_i$ that exits or enters into $\Gamma$ must pass through $A$).
The base space $B^d (\Gamma, A,n)$ is contractible since there is a deformation retraction of $\Gamma$ onto $A$ which extends to $B^d(\Gamma, A,n)$.

Next write an element in $B^d (\Gamma, A,n)$ as an equivalence class
$[[x_1,\ldots, x_k]]$ with $x_i\in\Gamma - A$ and some $k\leq n$. The preimage
$\lambda^{-1}[[x_1,\ldots, x_k]]$ consists of all possible unordered $n$-tuples containing
$x_1,x_2,\ldots , x_k$ with remaining entries $y_1,\ldots, y_{n-k}$ such that
$[y_1,\ldots, y_{n-k}]\in B^d((V-\Gamma)\cup A,n-k)$. This preimage is a copy of
$B^d((V-\Gamma)\cup A, n-k)$.
By induction on the number of simplices of $V$, we can assume that
$B^d((V-\Gamma)\cup A, n-k)$ is $2d-2$-connected (the case of a single simplex has been discussed at the beginning of the proof).
The map $\lambda : B^d (V,n)\lrar B^d (\Gamma, A,n)$ has then a contractible base and preimages that are $2d-2$-connected. We wish to show that the total space is $2d-2$-connected. We cannot use the Smale-Vietoris theorem (Theorem \ref{smale}) directly since $\lambda$ is not proper. To go around this, we pass to the compactified versions and show that $\lambda$ can be deformed to a proper map.
Let ${\mathcal W}^d(V,n)$ be the compact deformation retract of $B^d(X,n)$ as  discussed in Corollary \ref{mathcalw}. The restriction $\tilde\lambda$ of $\lambda$ to ${\mathcal W}^d(V,n)$ maps onto
${\mathcal W}^d(\Gamma,A,n)\subset B^d(\Gamma, A,n)$ and
we have the diagram
$$\xymatrix{
{\mathcal W}^d(V,n)\ar[d]^{\tilde\lambda}\ar[r]&B^d(V,n)\ar[d]^\lambda\\
{\mathcal W}^d(\Gamma, A,n)\ar[r]&B^d(\Gamma, A,n)
}$$
where the horizontal maps are inclusions and deformation retractions.
 This last statement follows from the fact that the deformation retraction of $B^d(V,n)$ onto ${\mathcal W}^d(V,n)$ descends to a deformation retraction of $B^d(\Gamma, A,n)$ onto ${\mathcal W}^d(\Gamma, A,n)$ as a consequence of Lemma \ref{property2}.
Thus given a configuration $\zeta= [[x_1,\ldots, x_k]]\in {\mathcal W}^d(\Gamma, A,n)$, $k\leq n$, $x_i\not\in A$, we can consider its preimage $\tilde\lambda^{-1}(\zeta )$ in
${\mathcal W}^d(V,n)$ and its preimage $\lambda^{-1}(\zeta )$ in $B^d(V,n)$. Then
$\lambda^{-1}(\zeta )$ deformation retracts onto $\tilde\lambda^{-1}(\zeta )$.
Here
$\tilde\lambda^{-1}(\zeta )\subset \lambda^{-1}(\zeta)= B^d((V-\Gamma)\cup A, n-k)$. Since $\lambda^{-1}(\zeta )$ is $2d-2$-connected, this shows that $\tilde\lambda^{-1}(\zeta )$ is also $2d-2$-connected. The map $\tilde\lambda$ is now proper being a map between compact spaces. Moreover both total and base spaces are
connected by Lemma \ref{connected}. We can invoke Theorem \ref{smale} to conclude that the total space
${\mathcal W}^d(V,n)$ and hence $B^d(V,n)$ are $2d-2$-connected as desired.
\end{proof}


\section{Proof of Theorem \ref{main2}}\label{proofthm2}

Our objective is to find conditions on $X$ so that the inclusion $\Delta^d(X,n)\hookrightarrow X^n$ induces an isomorphism on some homotopy groups through a range (the ``homotopical depth"). The proof given in the unordered case $B^d(X,n)$ (section \ref{proofofmain1}) fails here because the analogue of \eqref{lambdamap} is now a map
$\Delta^d (V,n)\lrar \Delta^d (\Gamma, A,n)$ which has \textit{disconnected} fibers so Smale's theorem doesn't automatically apply. In fact we need an entirely new approach.

First some definitions.

\bde\label{defs}\
\begin{itemize}
\item If $x\in X$, $U$ a neighborhood of $x$, then we call $V$ a sub-neighborhood (of $x$ in $U$) if $V$ open, $x\in V\subset U$.
\item A space $X$ is \textit{locally contractible} if for any $x\in X$ and any neighborhood $U$ of $x$, there is a sub-neighborhood $V$ which deformation retracts onto $x$.
\item A space $X$ has local homotopical dimension $k$ if $x,U$ as above, there is a subneighborhood $V$ with $V-\{x\}$ being $k$-connected. For instance being \textit{locally punctured connected} means having $0$ local homotopical dimension. A manifold of dimension $m$ has local homotopical dimension $m-2$ but not $m-1$.
\end{itemize}
\end{definition}

If $X$ is a simplicial complex, we call a ``chamber" of $X$ any simplex that is not contained in another simplex as a face. Obviously if $X$ has local homotopical dimension $r$, then chambers must have dimensions at least $r+2$. We call  a simplex a ``shared face" if it is shared by two chambers or more. This shared face doesn't need to be of codimension one. In figure \ref{contractnhbd}, the complex on the far right is made out of three chambers (of dimension $2$) joining along a shared edge. A shared face $A=\Gamma_1\cap\cdots\cap\Gamma_k$
is called essential if $X=\Gamma_1\cup\cdots\cup\Gamma_k$ is not a cell (i.e. homeomorphic to a ball or to a halfball). This rules out cases like $X$ being a regular polygon  triangulated so that the origin $A=o$ is the common vertex of all triangles. A neighborhood $V\setminus \{o\}$ is up to homotopy a circle in this case, so that $o$ behaves like an interior point of a chamber (and is inessential).

\ble\label{formadmiss}
A finite simplicial complex $X$ has local homotopical dimension $r$ if and only if all chambers are of dimension at least $r+2$ and all essential shared faces of dimension at least $r+1$.
\ele

\begin{proof}
It suffices to consider points $x\in X$ that are either in the interior of a chamber or in the interior of a shared face.
In the case $x$ is in the interior of a chamber, $V\cong\bbr^m$ so $m$ (the dimension of the chamber) must be at least $r+2$ (Proposition \ref{bdrmn}). On the other hand, if $x$ lies in the interior of a shared face $A$,
a small neighborhood $V$ of $x$ is the union of chambers $\Gamma_1\cup\cdots\cup\Gamma_q$ joining along $A$, with $q\geq 2$. If $A$ is inessential, then a neighborhood $V$ of $x\in A$ is either a ball or a half-ball of dimension at least $r+2$. Suppose $x$ essential and let $s=\dim A$. Then $V-\{x\}\simeq \bigvee S^{s}$ is a bouquet (this holds even if $s=0$ and $A$ is vertex). Since this neighborhood must be $r$-connected, $s$ must be at least $r+1$.
\end{proof}

The following is our main statement.
Here we assume $d<n$ since otherwise $\Delta^d(X,n)=X^n$ and there is nothing to prove.

\bth\label{depth}
Let $X$ be a locally finite polyhedral space with local homotopical dimension $r$, $r\geq 0$, and  let $1\leq d<n$. Then $$\pi_i(\Delta^d(X,n))\cong \pi_i(X)^n\ \ \ \hbox{for}\ \ i\leq rd + 2d-2\ .$$
\end{theorem}

\begin{proof} The starting point is Lemma \ref{key}. As in the proof of Lemma \ref{formadmiss}, a contractible neighborhood $V$ of $x\in X$ is of three types: (i) $V\cong\bbr^m$ with $m\geq r+2$, (ii)  half-space $H$ of dimension $m\geq r+2$ or (iii) it is a union of such half-spaces along a shared face of dimension at least $r+1$. We must show that $\Delta^d(V,n)$ is $dr+2d-2$ connected.

In the case $V\cong\bbr^m$ we know by Corollary \ref{manifold1} that $\Delta^d(\bbr^m,n)$ is $dm-2$-connected, and that $dm-2=d(r+2)-2\geq dr+2d-2$ as claimed.
If $V$ is homeomorphic to half-space $H$ (with boundary $L$), $V$ is the direct limit of a nested sequence of spaces $H_i\supset H$ such that $H_i\cong\bbr^m$ for all $i$, it follows that $\Delta^d(H,n)$ is the direct limit of spaces
$\Delta^d(H_i,n)$ that are $dr+2d-2$-connected, and hence has the same connectivity at least.
We are left with the case $V$ is the union of simplices (chambers) $\Gamma_1\cup\cdots\cup\Gamma_q$ joining along an essential face $A$. We can assume wlog that any two faces join along $A$, i.e.
$\Gamma_i\cap\Gamma_j=A$.
Luckily the structure of this neighborhood $V$ is sufficiently nice to allow us to give a decomposition of $\Delta^d(V,n)$ as the colimit of an explicit diagram.

We start by observing that each configuration of $n$ points of $V$ gives rise to a tuple of integers $(k_1,\ldots, k_q)$, $k_1+\cdots + k_q \geq n$, where $k_i$ denotes the number of points of the configuration inside the face $\Gamma_i$. Obviously these $k_i$-configurations can overlap when points of the configuration fall in $A$. Keeping track of the various overlaps can be expressed in terms of a poset of intersections. More precisely set the index set $$I = \{1,2,\ldots, q\}^n = \{(i_1,\ldots, i_n)\ |\ i_j\in \{1,2,\ldots, q\}\} .$$
We can cover $\Delta^d(V,n)$ by the closed sets
$U_{(i_1,\ldots, i_q)}$, $(i_1,\ldots, i_q)\in I$, where
$$U_{(i_1,\ldots, i_q)} =\{(x_1,\ldots, x_n)\in \Delta^d(V,n)\ |\  x_j\in \Gamma_{i_j} , \ i_j\in
\{1,2,\ldots, q\}\}.$$
Let $\mathcal D$ be the intersection poset $P_U$ associated
to the cover $U_I$ of $\Delta^d(V,n)$, also referred to as \textit{subspace diagram}. It is clear by construction that
colim$\mathcal D$ is precisely $\Delta^d(V,n)$. Here's how this poset diagram looks
like for $k=2$ and $d=1$, i.e. for the configuration space $F(X\cup_A Y, 2)$ (see \cite{sun}, Theorem 2.0.17):
\begin{equation}\label{diagram}
\xymatrix@C=-0.9em{
&X\times Y-\Delta A&&F(Y,2)&&Y\times X-\Delta A&&F(X,2)\\
X\times A-\Delta A\ar[ur]&&A\times Y-\Delta A\ar[ur]\ar[ul]&&
Y\times A-\Delta A\ar[ur]\ar[ul]&&
A\times X-\Delta A\ar[ur]\ar[ul]&&
X\times A-\Delta A\ar[ul]\\
&&&&F(A,2)\ar[urrrr]\ar[rru]\ar[u]\ar[ull]\ar[ullll]}
\end{equation}
(the spaces on the extreme right and left are being identified).

Going back to the general diagram $\mathcal D$, since all inclusions are closed cofibrations (this is standard to check \cite{sun}), we have
$$\Delta^d(V,n)=\hbox{colim}(\mathcal D)\simeq \hbox{hocolim}(\mathcal D) .$$
In fact the canonical map from the homotopy colimit of a sequence of inclusions of T1 topological spaces to the actual colimit is a weak equivalence (see \cite{kriz}).
The connectivity of this (sequential) hocolim
$\Delta^d(V,n)$ is at least the least connectivity of the spaces making
up the diagram. If we set $\Gamma_0=A$, these spaces are of
the form $\Gamma_{i_1}\times\cdots\times \Gamma_{i_n} \cap\Delta^d(V,n)$
(we refer to these subspaces as the ``constituent subspaces" of the diagram).
Each of these constituent subspaces
is quite manageable and we can apply the localization principle to it.
Indeed $\Gamma_{i_1}\times\cdots\times \Gamma_{i_n} \cap\Delta^d(V,n)$ is the complement in
$\Gamma_{i_1}\times\cdots\times \Gamma_{i_n}$ of subspaces of certain codimensions. The smallest codimension is attained by $\Delta_{d+1}(A,n)$ in $A^n$ (that is for $\Delta^d(A,n)$).
If $s=\dim A$, then this codimension is $ds$. It follows that the smallest connectivity among the constituent subspaces  is $ds-2\geq d(r+1)-2=dr+d-2$.
As pointed out the connectivity of $\Delta^d(V,n)$ (as a hocolim) must be at least the connectivity of $\Delta^{d}(A,n)$ which is $dr+d-2$.
This is not quite the connectivity we seek and we  must improve it by $d$.

To do so observe that there is associated
to the poset $P_U$ of the cover a natural filtration whose $j$-th space is
${\mathcal F}_j = colim P_{j}$ where $P_j$ is the poset consisting of
 \begin{eqnarray*}
\Gamma_{i_1}\times\cdots\times\Gamma_{i_n}\cap\Delta^d(V,n)
\ , & i_{k_1}=\cdots = i_{k_{s}} = 0\ \hbox{for}\ s\geq n-j\\
&\hbox{and some subset}\ \
\{k_1,\ldots, k_{s}\}\subset \{1,\ldots, n\}
\end{eqnarray*}
with $i_s\in \{0,1,\ldots, q\}$ and $\Gamma_0=A$ as pointed out.
In other words, $\mathcal F_j$ is the
subspace where at most $j$ of the entries can be outside $A$. We have the series of inclusions
$${\mathcal F}_0=\Delta^d(A,n)\subset {\mathcal F}_1\subset\cdots\subset {\mathcal F}_n=\Delta^d(V,n) .$$

 If we organize our poset vertically as in \eqref{diagram}, then
 $\mathcal F_j$ is the pushout of the first $j+1$ rows from the bottom.

For example $F(\Gamma_1\cup_A\Gamma_2,2)$ (the case depicted in diagram \eqref{diagram} with
$X=\Gamma_1,Y=\Gamma_2$), there are three filtration terms starting with
$\mathcal F_0=F(A,2)$, the colimit $\mathcal F_1$
of the first two rows and ${\mathcal F}_2$ being the whole colimit. The special case of $F(\bbr^2,2)=\Delta^1(\bbr^2,2)$ is enlightening ($d=1$, $n=2$), where here we write $\bbr^2= \Gamma_1\cup_A\Gamma_2$ with $\Gamma_i$'s being two halplanes joining along $A\cong\bbr$. The first filtration term is $\mathcal F_0 = F(\bbr,2)\simeq S^0$. The next filtration term is
$$\mathcal F_1= (\Gamma_1\times A\cup A\times\Gamma_1\cup
\Gamma_2\times A\cup A\times\Gamma_2)\cap \Delta^1 (\bbr^2,2). $$
Each term $(\Gamma_i\times A\cup A\times\Gamma_j)\cap \Delta^1 (\bbr^2,2)=\Gamma_i\times A\cup A\times\Gamma_i-$diag$(A)$ deformation retracts onto
a circle so $\mathcal F_1$ is the union of two circles along an $S^0$; i.e
$\mathcal F_1\simeq S^1\vee S^1\vee S^1$.
Finally $\mathcal F_2\simeq F(\bbr^2,2)\simeq S^1$. The connectivity has changed
going from $\mathcal F_0$ to $\mathcal F_1$. It has remained stable afterwards. This happens more generally; i.e. we will show that the connectivity jumps going from $\mathcal F_i$ to $\mathcal F_{i+1}$ for $0\leq i< d$.
More precisely in this range, the connectivity of $\mathcal F_k$ is $c_k=ds+k-2$.
We have already argued that the connectivity of
$\mathcal F_0=\Delta^d(A,n)$ is $ds-2$.

Let's organize into a row $R_k$ the constituent subspaces
$\Gamma_{i_1}\times\cdots\times\Gamma_{i_n}\cap\Delta^d(V,n)$ where precisely $k$ of the
$\Gamma_{i_j}$'s are not equal to $A=\Gamma_0$. One point we will capitalize on is that in the range $0\leq k\leq d$,
$\Gamma_{i_1}\times\cdots\times\Gamma_{i_n}\cap\Delta^d(V,n)$ is the complement in $\Gamma_{i_1}\times\cdots\times\Gamma_{i_n}$ of tuples with $d+1$-diagonal elements lying only in $A$.
At the first stage, all components of $R_1$ intersect
along $\mathcal F_0 = \Delta^d(A,n)$.

If $n=d+1$, the situation is very clear. Here
$\Delta^d(A,d+1)=A^{d+1}-$diag$(A)\simeq S^{ds-1}$, and all constituent subspaces
for $k\leq d$ are of the form
$\Gamma_{i_1}\times\cdots\times\Gamma_{i_{d+1}}\cap\Delta^d(V,d+1)=
\Gamma_{i_1}\times\cdots\times\Gamma_{i_{d+1}}-\hbox{diag}(A)$,
thus they are contractible since they are the complement of a closed subspace in the boundary
of a cube. This means that going up the filtration, we are suspending in various ways the
spherical class, as in the example discussed earlier, and the connectivity in homology
is going up by one at every step.

For more general $n$, the constituent subspaces are
not in general contractible but we have the useful lemma

\ble\label{null} The inclusion $\mathcal F_{k-1}\hookrightarrow \mathcal F_k$ is null-homotopic for $k\leq d+1$.
\ele

\begin{proof} We need some notation.
We introduce $\mathcal F_{k}(n)$ for the filtration terms of $\Delta^d(V,n)$ (we added the
index $n$ to the previous notation). We also introduce
 $\mathcal F_{k,j}(n)$ for the subspace of all configurations $(x_1,\ldots, x_n)\in\mathcal F_{k}(n)$ where $x_j$ can be in all of $H$.
 We have that $\mathcal F_k(n)=\bigcup_{1\leq j\leq n}\mathcal F_{k,j}(n)$.
There is an inclusion
$$\mathcal F_{k-1}(n)\hookrightarrow \mathcal F_{k,n}(n)\subset\mathcal F_k(n)$$
On the other hand there are various
emdeddings of
$\mathcal F_{j}(n)$ into $\mathcal F_{j}(n+1)$ one of which sends
\begin{equation}\label{embed}
(x_1,\ldots, x_{n})\longmapsto (\phi_1 (x_1),\ldots, \phi_1 (x_{n}), p_n)
\end{equation}
where $\phi_t$ is any isotopy of the halfspace $H$ extending an isotopy
of $A$ onto its halfspace $(a_1,\ldots, a_s)$, $a_1< 0$, and where
$p_n = (n,0,\ldots, 0)\in A\cong\bbr^s$.
The first observation is that the inclusion
$\mathcal F_{k-1}(n)\hookrightarrow \mathcal F_{k,n}(n)$
is homotopic to the composite
$$\mathcal F_{k-1}(n)\lrar
\mathcal F_{k-1}(n-1)
\hookrightarrow
\mathcal F_{k-1}(n)\hookrightarrow
\mathcal F_{k,n}(n)$$
where the first map is projection discarding the last configuration, and
where the middle map is the inclusion \eqref{embed}. The idea here is that the last
coordinate $x_n\in H$ can be moved in $H$ away from $A$, all configurations
are then mapped by $\phi_t$, and after that the last coordinate if brought down to
$p_n$. Note that the last configuration can move in $H$ without constraint since
$k\leq d$. Next we factor the composite above $\mathcal F_{k-1}(n-1)
\hookrightarrow
\mathcal F_{k}(n)$ through
$\mathcal F_{k-1}(n-1)
\hookrightarrow
\mathcal F_{k,n-1}(n)$ and reiterate this construction to factor the map
up to homotopy this time through $\mathcal F_{k-1}(n-2)$, etc. At the end, the map
$\mathcal F_{k-1}(n)\hookrightarrow \mathcal F_k(n)$ factors through
$\mathcal F_{k-1}(k-1)$ which is contractible.
\end{proof}

Lemma \ref{null} implies that for $k\leq d+1$ and up to homotopy
\begin{equation}\label{susp}
\mathcal F_k/\mathcal F_{k-1}\simeq \mathcal F_k\vee\Sigma\mathcal F_{k-1}.
\end{equation}
Set $I = (i_1,i_2,\ldots, i_n)$ an ordered tuple with $n-k$ entry $0$, and write
\begin{equation}\label{vi}
V_I:=\Gamma_{i_1}\times\cdots\times\Gamma_{i_{n}}\cap\Delta^{d}(V,n) ,
\end{equation}
which is a constituent subspace for the $k$-th row.
Write $V_{\hat{I}}=V_I\cap \mathcal F_{k-1}$. This consists of tuples
obtained from $V_I$ by replacing a non-zero entry of $I$ by $0$. Then
$\mathcal F_k/\mathcal F_{k-1}\simeq\bigvee_I V_I/V_{\hat{I}}.$ If we show
that for any $I$, $V_I/V_{\hat{I}}$ is one degree more connected than $V_I$, then necessarily
$\mathcal F_k$ gains one degree of connectivity by \eqref{susp}. The argument here is
similar to \cite{dt}, \S4 and \S5. Assume wlog that $V_I=\Gamma_1\times\cdots\times\Gamma_n\cap\Delta^d(V,n)$ with $\Gamma_n=H$, and $\Gamma_i=A$ or $H$ for $i<n$. Project onto the first $n-1$ factors
$V_{I}\lrar \Gamma_1\times\cdots\times\Gamma_{n-1}\cap\Delta^d(V,n-1)$.
Call this target space $B_I$. It is $c_{k-1}=ds+k-3$ connected by induction. A chain in $V_I$ can be seen
as a chain in $B_I\times H$ transversal to some forbidden subspaces, and this is a
cone on a chain in $B_I$. This cone becomes a suspension class in $V_I/V_{\hat I}$.
If we write h-connected for homological connectedness, then
the h-connectivity keeps increasing by $1$ at every stage, for
$k<d$ so that $\mathcal F_k/\mathcal F_{k-1}$, and thus $\mathcal F_k$, is $ds+k-2$ connected . Since $s\geq r+1$, this is $dr+2d-2$ h-connected. This connectivity will not go down going up the filtration.

Finally to get the connectivity, we need argue that $\mathcal F_n$ is simply connected.
In fact $\mathcal F_k$ becomes $1$-connected as soon as
$k\geq 1$. To see this, we go back to the colimit diagram \eqref{diagram} where the smallest connectivity of the constituent subspaces $V_I$ as in \eqref{vi} is $d(r+1)-2$. When this is larger than $1$, each $V_I$ is
simply connected, so is the colimit and the theorem holds.
Now some $V_I$ fail to be simply connected when
$d(r+1)\leq 2$, that is when: (i) $r=0,d=1$, (ii) $r=1=d$, and (iii) $r=0,d=2$. In the first
case, the theorem is equivalent to saying that $F(X,n)$ is connected if $X$ is locally punctured connected. This is precisely Lemma \ref{connected} so this case is settled. In the case (ii) we are looking at $\Delta^d(A,n)=\Delta^1(\bbr^2,n)=F(\bbr^2,n)$ as the bottom space of our colimit
diagram. This is of course not simply connected, but the map $\pi_1(\mathcal F_0)\rightarrow\pi_1(\mathcal F_1)$ is the trivial map since induced from a null homotopic map (Lemma \ref{null}), so that $\mathcal F_1$, and inductively $\mathcal F_k$, are simply connected by Van-Kampen. Remains the case (iii) which occurs when $\Delta^d(A,n)=\Delta^2(\bbr,n)$. The fundamental group of this space is discussed in Example \ref{khovanoex}. Here too the fundamental group trivializes from $\mathcal F_1$ onwards so that $\mathcal F_n=\Delta^d(V,n)$ is simply connected.
\end{proof}

\subsection{The Homology of the Filtration Terms}
This subsection is of independent interest and
gives a useful description of the homology of the filtration terms. This is
sketchy but details can be filled in.
First of all there is a nice way to see that the inclusion
$\mathcal F_0\rightarrow\mathcal F_1$ induces the trivial map in homology without resorting to Lemma \ref{null}.
Here $\mathcal F_0=\Delta^d(A,n)$ has torsion free homology admitting a basis realized by products of spheres \cite{dt}. We need understand how these homology classes occur. There is a spherical class in $\Delta^d(A,d+1)=A^{d+1}-$diag$(A)\simeq S^{ds-1}$, where $s=\dim A$. Now
$\Delta^d(A,d+1)$ embeds in $\Delta^d(A,n)$ in many ways as in \eqref{embed} (recall that $d<n$).  This embedding has a retract so induces a monomorphism in homology. The image of the spherical class in this case is denoted $\{x_1,\ldots, x_{d+1}\}$. The various other embeddings obtained by choosing another subset $\{i_1,\ldots, i_{d+1}\}\subset\{1,2,\ldots, n\}$ where to insert the $a$'s give rise to spherical homology classes
$\{x_{i_1},\ldots, x_{i_{d+1}}\}$. These classes generate the homology of $\Delta^d(A,n)$ in a very precise sense. There is an action of the operad $\{D^s(k)\}_{k\geq 0}$ of little $s$-dimensional disks on
$\bigcup_{n\geq 1}\Delta^d(A,n)$, where $s=\dim A$ and $D^s(k)$ is the space of $k$ pairwise disjoint open disks in the unit disk of dimension $s$ (to keep with the terminology the word ``disk" is used instead of ``ball"). The action of $D^s(2)\simeq S^{s-1}$ is given as follows
$$D^s(2)\times \Delta^d(A,n_1)\times \Delta^d(A,n_2)\lrar \Delta^d(A,n_1+n_2)$$
and yields a bracket operation in homology
$$[-,-]: H_p(\Delta^d(A,n_1))\otimes H_q(\Delta^d(A,n_2))\lrar
H_{p+q+s-1}(\Delta^d(A,n_1+n_2))$$
The product map in homology is given by the action of $H_0(D^s(2))$ and is the induced
map in homology of the concatenation of two configurations after placing the first one
in a disk of radius $1/2$ centered at $(-1/2,0,\ldots, 0)$ and the other in another
disk of same radius centered at $(1/2,0,\ldots, 0)$.
One main theorem of \cite{dt} reads as follows. The bracket of two cycles
is important to understand and can be described as follows. Given a cycle (or a chain) $c$
in $\Delta^d(A,n)$, we say we ``localize it" in a disk $D^s$ if we choose a homeomorphism (which
can be made canonical) between $A\cong\bbr^s$ and $D$, and take the image of $c$ in
$\Delta^d(D,n)$ via this homeomorphism. We obtain the bracket $[\alpha_1,\alpha_2]$ by localizing
the cycles respectively in two disjoints disks $D_1$ and $D_2$ and taking the new cycle
obtained by rotating $D_1$ around $D_2$ (or $D_2$ around $D_1$, up to sign) in $\bbr^s$.

\bth\label{victor} (Proposition 3.9 of \cite{dt})
The homology of $\Delta^d(A,n)$ is torsion free, generated additively by products of iterated
brackets where each factor is either $x_i$ or an iterated
bracket of the form
$$[\ldots [[B_1,B_2],B_3]\ldots, B_\ell]\ ,\ \ell\geq 1$$
where each $B_s$ is of the form
$$B_s = [\ldots [[\{x_{j_{1,s}}, x_{j_{2,s}},\ldots, x_{j_{d+1,s}}\}, x_{i_{1,s}}], x_{i_{2,s}}]\ldots, x_{i_{\ell_s,s}}]$$
(further conditions are stated on indices to get a basis).
\end{theorem}

Let's argue for example that
$\{x_1,\ldots, x_{d+1}\}$ maps to zero in the homology of the next filtration term. Consider the diagram of inclusions
$$\xymatrix{
\Delta^d(A,d+1)\ar[r]\ar[d]^\iota&\Delta^d(A,n)=\mathcal F_0\ar[d]\\
(\Gamma_1\times A^d)\cap \Delta^d(V,n)\ar[r]&\mathcal F_1
}$$
The bottom space $(\Gamma_1\times A^d)\cap \Delta^d(V,n)=\Gamma_1\times A^d-$diag$(A)$ is contractible since we are removing a subspace from the boundary of $\Gamma_1\times A^d$.
The map $\iota$ is trivial and the commutativity of the diagram shows that
$\{x_1,\ldots, x_{d+1}\}$ maps trivially in $\mathcal F_1$. A class of the form
$[\{x_1,\ldots, x_{d+1}\}, x_{d+2}]$ dies in $\mathcal F_1$ for example since this class can
be represented by the composite
\begin{equation}\label{threemaps}
S^{ds-1}\times S^{d-1}\lrar A^{d+2}-\hbox{sing}\hookrightarrow
A^{d+1}\times H-\hbox{sing}\hookrightarrow\mathcal F_1
\end{equation}
The first map is obtained from the operadic action. Here the factor $S^{d-1}$ is the locus of
$x_{d+2}$ rotating in some sphere in $\bbr^d$, so when $x_{d+2}$
is  allowed to be in $H$, this sphere is coned off and the composite of the first two maps
in \eqref{threemaps} is trivial on the top homology class which by definition is
$[\{x_1,\ldots, x_{d+1}\}, x_{d+2}]$. A similar argument applies to show that the image of
$B_s$ as in the notation of Theorem \ref{victor} is trivial in $\mathcal F_1$.
For the image of the bracket $[B_s,B_t]$ one can argue similarly. One constructs this class by
localizing $B_s$ and $B_t$ in distinct disks $D_1$ and $D_2$, and rotating one disk around the other. But the class $B_s$ is the boundary of a chain in
$H'\times D_1^{n-1}\cup D_1\times H'\times D^{n-2}\cup\cdots\cup D^{n-1}\times H'\subset\mathcal F_1$,
where $H'$ is the part of $H$ with boundary $D_1$.
This means that $[B_s,B_t]$ must map to zero in $H_*(\mathcal F_1)$.
Remains to see that the image of a product is trivial but this is immediate.

Note that there are many ways
a given class $[\ldots [[B_1,B_2],B_3]\ldots, B_\ell]$
can die in $\mathcal F_1$, and so
in $\mathcal F_1$ we obtain suspension classes one degree higher. This describes
completely the homology of $\mathcal F_1$ and obviously it is one degree more
connected than $\mathcal F_0$.


\section{Fundamental Groups}\label{fundamental}

 In this final section we take a more pedestrian look at the isomorphism $\pi_1(B^d(X,n))\cong H_1(X,\bbz )$ for $d\geq 2$. This is expressed in terms of braids.
 As before $X$ is a simplicial complex. Note that loops in $\sp{n}X$,  based at a basepoint of the form $[\ast,\ldots,\ast]$ say, lift to $X^n$ under the quotient projection (see \cite{kt}, section 5 for example)
$$\xymatrix{&X^n\ar[d]\\
S^1\ar[r]\ar@{..>}[ru]&\sp{n}X
}$$
This says that a homotopy class of a loop
$\gamma : S^1\lrar\sp{n}X$  based at $[\ast,\ldots,\ast]$ can be represented by a tuple $[\gamma_1,\ldots,\gamma_n]$,
where $\gamma_i:S^1\lrar X$ is a loop in $X$. Moreover and by the simplicial approximation theorem, any loop in $X^n$ deforms into an $n$-tuple of simplicial loops in $X$ so that $\gamma$ can be represented by an unordered tuple of simplicial loops in $\sp{n}(X)$ for some simplicial decomposition.

We can try to describe loops in $\Delta^d(X,n)$ and $B^d(X,n)$ in the same way but both spaces are \text{not} simplicial complexes in general, only of the homotopy type of one. However after passing to a barycentric subdivision, $B^d(X,n) = \sp{n}X-B_{d+1}(X,n)$ deformation retracts onto a cellular complex ${\mathcal W}^d(X,n)$ (Lemma \ref{property1}). A loop $S^1\lrar B^d(X,n)$ deforms into a loop into ${\mathcal W}^d(X,n)$ which is cellular. Therefore and without loss of generality, we can represent a loop $\gamma : S^1\lrar B^d(X,n)$ within its homotopy class by a tuple of paths $t\longmapsto [\gamma_1(t),\ldots, \gamma_n(t)]$, with $\gamma_i$ a simplicial path in $X$ (not necessarily a closed loop) and $t\in [0,1]$. This is a \textit{braid with $n$-strands}. These paths or strands at any time $t$ do not intersect in more than $d$ points, and $[\gamma_1(0),\ldots, \gamma_n(0)] = [\gamma_1(1),\ldots, \gamma_n(1)]$. Similarly for loops into $\Delta^d(X,n)$.

\begin{figure}[htb]
\begin{center}
\epsfig{file=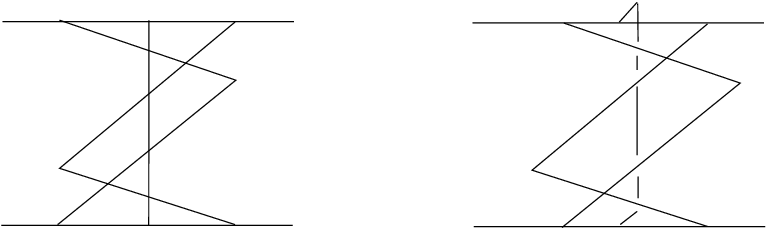,height=0.8in,width=2in,angle=0.0}
\caption{}\label{nullhomotopy}
\end{center}
\end{figure}

As a first example, consider $X=I$ the unit interval. By codimension argument, $\Delta^d(I,n)$ is simply connected if $d\geq 3$, so the only interesting case is when $d=2$ and we are removing from $I^n$ codimension $2$ subspaces corresponding to when $x_i=x_j=x_k$.
According to Example \ref{casestudy}, $\Delta^2(I,3)\simeq S^1$ and $\pi_1(\Delta^2(I,3))\cong\bbz$. An element in the fundamental group can be represented by a braid with $3$-strands embedded in $I\times I$ not all of which can pass by the same point at the same time. A non-trivial element is depicted in the lefthand side of figure \ref{nullhomotopy}. This braid cannot be trivialized in $\Delta^2(I,3)$ and it is amusing to try to test this fact. By moving the strands around keeping their endpoints fixed, there is no way we can separate them without going through a triple point.

\bex\label{khovanoex}
For $n\geq 3$, the fundamental group of $\Delta^2(I,n)$ has been analyzed by Khovanov \cite{khovanov}. There he shows that $\Delta^2(I,n)$ is a $K(\pi,1)$ and then gives a presentation for $\pi$. This presentation is given as follows. Define the right-angled Coxeter group $TW_n$ to be the group generated by the simple transpositions $s_i = (i,i+1), i\in [n-1]$, subject to the relations
$$s_i^2 = 1\ \ \ ,\ \ \ s_is_j = s_js_i\ \hbox{if}\ |i-j|>1\ .$$
Define the map $\phi: TW_n\lrar \sn$ by $\phi(s_i)=s_i$ for all $i\in [n-1]$. Then
$\pi_1(\Delta^2(I,n))\cong\hbox{ker}\ \phi$.
\eex

\begin{figure}[htb]
\begin{center}
\epsfig{file=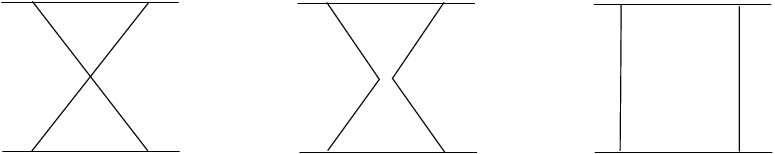,height=0.7in,width=2.5in,angle=0.0}
\caption{Resolving the intersection points}\label{resolve}
\end{center}
\end{figure}

In the unordered case it is possible to kill the ``braiding'' by \textit{interchanging} strands. Represent an element of $\pi_1(B^d(X,n))$ by a braid with $n$-strands embedded in $X\times I$. Suppose we have two intersecting strands. There is a way of ``resolving'' the intersection points as illustrated in Figure \ref{resolve}. The figure depicts a loop $f(t) = [f_1(t), f_2(t)]$ with two strands crossing for some $s\in [0,1]$. Define $\tilde f = [\tilde f_1,\tilde f_2]$ to be such that $\tilde f_i(t)=f_i(t)$ if $t\leq s$, and $\tilde f_1(t) = f_2(t)$,
$\tilde f_2(t)=f_1(t)$ if $t\geq s$. These give two representations of the same loop in $\Omega B^d(X,n)$ for $d\geq 2$. The difference however is that after changing $f$ by $\tilde f$, by a small homotopy we can now separate the strands of $\tilde f$ so that no intersection occurs. This also explains why the fundamental group must be abelian (\cite{kt})

For example using this resolution of intersections, we can immediately trivialize the braid in $\Omega B^2(\bbr,3)$ depicted on the left of figure \ref{nullhomotopy}. This is no surprise since $B^2(\bbr,3)\simeq B^2(I,3)$ is contractible and is identified with the three simplex with one edge removed.

The resolution of intersections when applied to loops in $\Omega V$, with $V$ a tree, implies that we have a surjection $\pi_1(B(V,n))\lrar\pi_1(B^d(V,n))$. Since $B^d(V,n)$ is connected for $d\geq 2$ (Lemma \ref{connected}), pick the basepoint in this fundamental group to be $[x_1,\ldots, x_n]$ with $x_i\neq x_j, i\neq j$, and write a braid $\gamma (t) = [\gamma_1(t),\ldots, \gamma_n(t)]$. As discussed, we can assume the $\gamma_i$'s to be \textit{non intersecting} strands. Since $V$ is one dimensional, then necessarily $\gamma_i(0)=\gamma_i(1) (=x_i)$ so all strands must start and finish at the same point. Each $\gamma_i$ can be homotoped to the constant strand at $x_i$, without further intersections, and the loop we started out with is trivial up to homotopy. The above discussion allows us to give a streamlined proof of the following proposition which we have already obtained as a Corollary to Theorem \ref{main1}.

\bpr If $X$ is a connected simplicial complex which is not reduced to a point, $n\geq 2, d\geq 2$, then there is an isomorphism $\pi_1(B^d(X,n))\cong H_1(X;\bbz )$.
\epr

\begin{proof}
We need show that the inclusion $B^d(X,n)\hookrightarrow\sp{n}X$ induces an isomorphism on fundamental group if $d>1$. If we invoke Lemma \ref{key} as before, this boils down to showing that for $V$ a contractible neighborhood in $X$, $B^d(V,n)$ is simply connected whenever $d\geq 2$ and for any $n\geq 1$.
If $V$ is a contractible neighborhood in simplicial $X$ (as in the proof of Theorem \ref{most}), any element in $\pi_1(B^d(V,n))$ can be represented by a braid and by resolving the intersection points, this braid can be homotoped to the trivial braid.
\end{proof}

\vskip 20pt

\bibliographystyle{plain[8pt]}

\end{document}